\documentclass[11pt,reqno]{amsart}

\usepackage{hyphenat}
\usepackage[a4paper,margin=31mm]{geometry}
\usepackage[foot]{amsaddr}
\usepackage{amsmath,amssymb,amsthm,mathtools}
\usepackage{microtype}
\usepackage{enumitem}
\usepackage[
  backend=biber,
  style=alphabetic,
  sorting=nyt,
  sortcites=true,
  maxbibnames=99,
  maxsortnames=99,
  maxalphanames=4,
  minalphanames=3,
  isbn=false
]{biblatex}

\AtBeginBibliography{\raggedright\interlinepenalty=10000}

\defbibheading{bibliography}[\refname]{\section*{#1}}
\usepackage{graphicx}
\usepackage{xcolor}
\usepackage[hypertexnames=false]{hyperref}
\usepackage{bookmark}
\usepackage{aliascnt}
\mathtoolsset{showonlyrefs}
\numberwithin{equation}{section}
\allowdisplaybreaks
\usepackage[section]{placeins}
\usepackage{tikz,booktabs}
\usetikzlibrary{arrows.meta}

\tikzset{
  tree vertex/.style={circle,draw=black,fill=white,line width=.4pt,
    minimum size=6.4mm,inner sep=0pt,font=\footnotesize},
  tree selected/.style={tree vertex,fill=black!12},
  tree edge/.style={draw=black,line width=.4pt},
  tree arrow/.style={-{Stealth[length=1.7mm,width=1.2mm]},
    draw=black,line width=.4pt},
  tree note/.style={font=\footnotesize,inner sep=1pt},
  tree ray/.style={tree arrow,draw=black!60,densely dashed}
}

\hypersetup{
  bookmarksdepth=2,
  colorlinks=true,
  linkcolor=blue!55!black,
  citecolor=green!45!black,
  urlcolor=blue!65!black,
  pdftitle={Height and diameter of split trees with an application to Plancherel trees},
  pdfauthor={Shengjun Zhang},
  pdfsubject={Height and diameter of split trees and Plancherel trees},
  pdfkeywords={Split trees, Plancherel measure, recursive Ewens trees,
    Poisson--Dirichlet distribution, branching random walks, height, diameter}
}

\theoremstyle{plain}
\newtheorem{theorem}{Theorem}[section]
\newaliascnt{proposition}{theorem}
\newtheorem{proposition}[proposition]{Proposition}
\aliascntresetthe{proposition}
\newaliascnt{lemma}{theorem}
\newtheorem{lemma}[lemma]{Lemma}
\aliascntresetthe{lemma}
\newaliascnt{corollary}{theorem}
\newtheorem{corollary}[corollary]{Corollary}
\aliascntresetthe{corollary}
\theoremstyle{definition}
\newaliascnt{remark}{theorem}
\newtheorem{remark}[remark]{Remark}
\aliascntresetthe{remark}
\newaliascnt{example}{theorem}
\newtheorem{example}[example]{Example}
\aliascntresetthe{example}
\theoremstyle{plain}

\newcommand{\Pp}{\mathbb P}
\newcommand{\E}{\mathbb E}
\newcommand{\T}{\mathcal T}
\newcommand{\U}{\mathcal U}
\newcommand{\cE}{\mathcal E}
\DeclareMathOperator{\Aut}{Aut}
\DeclareMathOperator{\GEM}{GEM}
\DeclareMathOperator{\PD}{PD}
\DeclareMathOperator{\Beta}{Beta}
\DeclareMathOperator{\Bin}{Bin}
\DeclareMathOperator{\Ewens}{Ewens}
\DeclareMathOperator{\Var}{Var}

\title[Height and diameter of split trees]
{Height and diameter of split trees
with an application to Plancherel trees}
\address{Université Paris-Saclay, Faculté des Sciences d'Orsay,
Institut de mathématiques d'Orsay, Bâtiment 307, F-91405 Orsay, France}
\email{zhang.shengjun@universite-paris-saclay.fr}
\author{Shengjun Zhang}

\begin{document}
\raggedbottom

\begin{abstract}
We study the height \(H_n\) and diameter \(D_n\) of a split tree after
the insertion of \(n\) balls. Each ball is stored at the first empty
vertex it reaches. The split vectors at the vertices are independent
copies of \(P=(P_i)_{i\ge1}\), which may have infinitely many positive
coordinates. Let \(M=\max_iP_i\). If
\(\Pp(M<1)>0\), then \(H_n/\log n\to c_\star\) and
\(D_n/\log n\to2c_\star\) almost surely and in \(L^p\) for every
\(1\le p<\infty\), where
\(c_\star=\inf_{\beta>1}\beta/(-\log\E\sum_i P_i^\beta)\).
If the infimum is attained at \(\beta_\star\), let
\(b_n=c_\star\log n-\tfrac{3c_\star}{2\beta_\star}\log\log n\).
Then \(H_n-b_n\) and \(D_n-2b_n\) are tight and satisfy uniform exponential
upper-tail bounds. For each \(p>0\), their \(p\)th absolute moments are
uniformly bounded if and only if \(\E[(-\log M)^p]<\infty\).
If \(\E[(-\log M)^p]\) is infinite, both moments are of order
\(\E[\min\{-\log M,\log n\}^{p}]\).
Each lower tail satisfies a uniform exponential bound if and only if
\(\E[M^{-\eta}]<\infty\) for some \(\eta>0\).
We apply these results to Fulman's Plancherel trees by taking
\(P\sim\PD(0,2)\).
\end{abstract}

\maketitle
\tableofcontents

\section{Introduction}\label{sec:introduction}

We study the height \(H_n\) and diameter \(D_n\) of a split tree with
\(n\) balls and one ball per occupied vertex. If a split assigns mass to more than one child
with positive probability, both height and diameter grow logarithmically.
When the infimum defining the height constant is attained, we show that
\(H_n-b_n\) and \(D_n-2b_n\) are tight for a centering sequence \(b_n\)
with a \(\log\log n\) correction. The mean height is \(b_n+O(1)\) if and
only if \(\E[-\log M]<\infty\), where \(M\) is the largest split
fraction. Remark~\ref{rem:mean-height} gives a split law for which
\(b_n-\E H_n\asymp\sqrt{\log n}\).

Our motivating example is a Plancherel measure on rooted trees.
The classical Plancherel measure assigns weight \((f^\lambda)^2/n!\) to a
partition \(\lambda\) of \(n\), where \(f^\lambda\) is the number of
standard Young tableaux of shape \(\lambda\). Under the
Robinson--Schensted correspondence, \(\lambda_1\) has the law of the
length of a longest increasing subsequence of a uniform random
permutation of \(\{1,\dots,n\}\). Vershik and Kerov~\cite{VershikKerov1977,VershikKerov1985} proved that
\(\lambda_1/\sqrt n\to2\) in probability. Logan and
Shepp~\cite{LoganShepp1977} independently showed that the limiting
constant is at least \(2\).

\begin{figure}[htbp]
\centering
\includegraphics[pagebox=cropbox,width=0.8\textwidth]{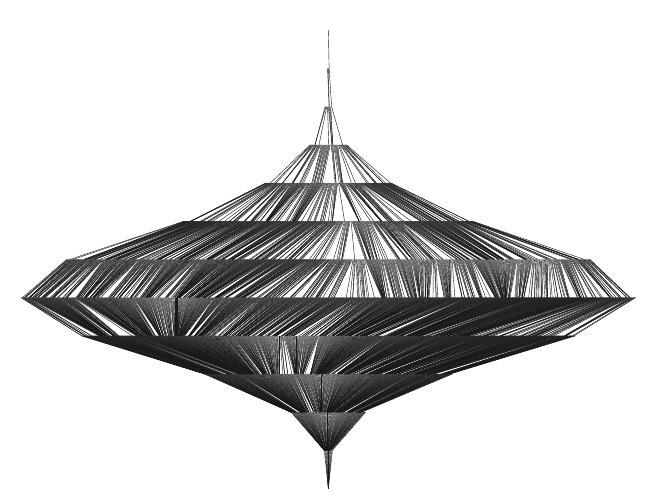}
\caption{A Plancherel tree with \(3000\) vertices and height \(12\).
The root is at the bottom.}
\label{fig:plancherel-tree}
\end{figure}

Fulman~\cite{Fulman2009} defined an analogue of this measure on rooted
trees, based on the growth and pruning operators of
Hoffman~\cite{Hoffman2003}, and proposed the study of its statistical
properties. In the spirit of Ulam's problem, we ask for the typical
height of a tree sampled from this measure.
A sample is shown in Figure~\ref{fig:plancherel-tree}.\footnote{Code for
generating and plotting Plancherel trees, written by Pierre-Lo\"ic
M\'eliot, is available at \url{https://github.com/darkxan771/trees}.}

We show that the height of such a tree with \(n\) vertices is
\[
 c_\star\log n-\frac{3c_\star}{2\beta_\star}\log\log n+O_{\Pp}(1),
\]
where \(\beta_\star\) and \(c_\star\) are given
by~\eqref{eq:plancherel-critical-equation} and \(O_{\Pp}(1)\) denotes
a term bounded in probability as \(n\to\infty\). The starting point is the
recursive description recalled in Theorem~\ref{thm:identification}.
The root subtrees correspond to the blocks of an Ewens partition of
\(n-1\) elements with parameter \(2\). Given the partition, they are
independent Plancherel trees, with \(|B|\) vertices in the subtree
corresponding to block \(B\).
The first version of this paper established this recursive description and proved the first-order height asymptotics in probability for the recursive Ewens trees. With independent Poisson--Dirichlet \(\PD(0,2)\) splits, the
split-tree construction of Devroye~\cite{Devroye1998} and
Janson~\cite{Janson2019} gives Fulman's upward
chain~\cite[Proposition~6.7]{Zhang2026}.
In Proposition~\ref{prop:plancherel-growth}, we compute the transition
probabilities directly and show that the Kuba--Panholzer insertion
procedure yields the same chain of tree shapes.

\subsection{Plancherel trees}\label{sec:plancherel-trees}

Let \(\T_n\) be the set of rooted unordered trees with \(n\) vertices,
up to root-preserving isomorphism. For a tree \(t\), we write \(|t|\) for its number of
vertices and \(t_v\) for the subtree formed by a vertex \(v\) and its
descendants. An \emph{increasing labelling} of \(t\in\T_n\) is a
bijection from its vertices to \(\{1,\dots,n\}\) that increases along
every path from the root. Let \(d(t)\) be the number of increasing
labellings of a representative of \(t\). Two labellings give isomorphic labelled trees if and only if they differ
by an automorphism of \(t\). Since the labels are distinct, only the
identity fixes a labelling. Thus each \(\Aut(t)\)-orbit has size
\(|\Aut(t)|\), and the number of increasingly labelled trees of shape
\(t\) is \(u(t)=d(t)/|\Aut(t)|\). The hook-length formula for
trees~\cite[Section~5.1.4, Exercise~20]{Knuth1973} gives
\begin{equation}\label{eq:hook}
 d(t)=\frac{|t|!}{\prod_{v\in t}|t_v|}.
\end{equation}
Fulman's \emph{Plancherel measure}~\cite[Definition~1]{Fulman2009} on
\(\T_n\) is
\begin{equation}\label{eq:plancherel}
 p_n(t)=\frac{d(t)\,u(t)}{Z_n},
 \qquad
 Z_n=\prod_{k=2}^{n}\binom k2=\frac{n!\,(n-1)!}{2^{n-1}}.
\end{equation}
In~\cite{Fulman2009} the numbers \(d(t)\) and \(u(t)\) are denoted by
\(m(t)\) and \(n(t)\). The fact that \(p_n\) is a probability measure is
due to Hoffman~\cite[Proposition~2.8]{Hoffman2003}, and a bijective
proof is recalled in Section~\ref{sec:double-labelling}. We write
\(T_n^{\mathrm{Pl}}\) for a random tree with law \(p_n\) and call it a
\emph{Plancherel tree}.

Fulman also defines an upward Markov chain \((X_n)_{n\ge1}\) whose
marginal at time \(n\) is \(p_n\)~\cite[Definition~2 and
Lemma~3.1]{Fulman2009}. Its transition probabilities are recalled in
Section~\ref{sec:pd-splits}, where the chain is identified with a
split-tree growth process. The \emph{height} \(H(t)\) of a tree \(t\)
is the largest distance from its root to a vertex, and its \emph{diameter} \(D(t)\) is the largest distance between two vertices.

\subsection{Recursive Ewens trees}\label{sec:ewens-trees}

For \(\theta>0\), the Ewens law \(\Ewens(m,\theta)\) assigns to a
partition \(\pi\) of a set with \(m\) elements the probability
\begin{equation}\label{eq:ewens}
 \frac{\theta^{\#\pi}}{\theta^{\overline m}}
       \prod_{B\in\pi}(|B|-1)!,
 \qquad
 \theta^{\overline m}=\theta(\theta+1)\cdots(\theta+m-1),
\end{equation}
where \(\#\pi\) is the number of blocks of \(\pi\) and
\(\theta^{\overline0}=1\)~\cite{Ewens1972,Pitman2006}. A
\emph{recursive Ewens tree} with parameter \(\theta\) on a finite set of
labels is constructed as follows. The root carries the smallest label. We partition the remaining labels
according to the Ewens law with parameter \(\theta\). Given this partition,
we repeat the construction independently on each block and attach the
resulting trees to the root. We write \(T_n^{(\theta)}\) for its shape
when there are \(n\) labels. For \(\theta=1\) the labelled tree is the
uniform random recursive tree (Section~\ref{sec:ewens-representation}).

\begin{theorem}[Plancherel trees are recursive Ewens trees]\label{thm:identification}
For every \(n\ge1\), the trees \(T_n^{\mathrm{Pl}}\) and \(T_n^{(2)}\)
have the same law. Equivalently, \(T_n^{\mathrm{Pl}}\) can be
constructed by sampling \(\pi\sim\Ewens(n-1,2)\) and attaching to a new
root, conditionally on \(\pi\), independent trees with laws
\(p_{|B|}\), \(B\in\pi\).
\end{theorem}

The proof is a short computation with~\eqref{eq:hook}
and~\eqref{eq:ewens}, given in Section~\ref{sec:ewens-representation}.
The recursive description is also implicit in the first proof
of~\cite[Theorem~4]{KubaPanholzer2012}, which counts double increasingly
labelled trees by a recursion on the sizes of the root subtrees.

\subsection{Height and diameter of split trees}\label{sec:height-intro}

Let \(P=(P_i)_{i\ge1}\) be a random probability vector, that is,
\(P_i\ge0\) and \(\sum_iP_i=1\) almost surely. We call \(P\) the
\emph{split vector}. To construct the split tree~\cite{Devroye1998,Janson2019},
give each vertex \(u\) of the infinite Ulam tree \(\U=\bigcup_{h\ge0}\{1,2,\dots\}^h\) an independent copy
\(P(u)\) of \(P\). Balls \(1,2,\dots\) are inserted one at a time at
the root. A ball that reaches an empty vertex is stored there, and a
ball that reaches an occupied vertex \(u\) moves to its \(i\)th child
with probability \(P_i(u)\). Let \(\widehat T_n\) be the tree of occupied
vertices after the insertion of \(n\) balls, and let
\(H_n=H(\widehat T_n)\), \(D_n=D(\widehat T_n)\) and \(M=\max_iP_i\).
We call the family of split vectors the
\emph{environment} and write \(\cE\) for the \(\sigma\)-field that it
generates. A formal construction is given in
Section~\ref{sec:coupling}. For \(\beta>0\) let
\[
 L(\beta)=\E\sum_{i\ge1}P_i^\beta\in(0,\infty],
 \qquad
 \kappa(\beta)=\log L(\beta),
\]
where \(\log\) denotes the natural logarithm.
If \(\Pp(M<1)>0\), then \(0<L(\beta)<1\) for every \(\beta>1\),
and we set
\begin{equation}\label{eq:transform-definition}
 c_\star=\inf_{\beta>1}\frac{\beta}{-\kappa(\beta)}.
\end{equation}
Lemma~\ref{lem:transform} shows that \(0<c_\star<\infty\).

\begin{theorem}[First-order height and diameter]\label{thm:height}
Assume that \(\Pp(M<1)>0\). Then, as \(n\to\infty\),
\[
 \frac{H_n}{\log n}\longrightarrow c_\star,
 \qquad
 \frac{D_n}{\log n}\longrightarrow 2c_\star
 \qquad\text{almost surely and in \(L^p\) for every \(1\le p<\infty\).}
\]
Moreover, almost surely, for every \(1\le p<\infty\),
\[
 \E\!\left[\left|\frac{H_n}{\log n}-c_\star\right|^p+
             \left|\frac{D_n}{\log n}-2c_\star\right|^p
             \,\middle|\,\cE\right]\longrightarrow0.
\]
\end{theorem}

If \(P\) is almost surely a unit mass, all balls follow the same ray and
\(H_n=D_n=n-1\). Now suppose that \(\Pp(M<1)>0\). The infimum
in~\eqref{eq:transform-definition} need not be attained.
Let \(p_{\max}=\|M\|_\infty\) be the essential
supremum of the largest coordinate. Lemma~\ref{lem:transform} shows
that \(L(\beta)^{1/\beta}\to p_{\max}\) as \(\beta\to\infty\), and that
\begin{equation}\label{eq:attainment}
 \text{the infimum in~\eqref{eq:transform-definition} is attained}
 \quad\Longleftrightarrow\quad
 \E\,\#\{i\ge1:P_i=p_{\max}\}<1 .
\end{equation}
In this case it is attained at a unique point
\(\beta_\star\in(1,\infty)\), and
\(c_\star=-1/\kappa'(\beta_\star)\). In particular, the infimum is
attained when \(p_{\max}=1\), as observed by
Joseph~\cite[discussion following Theorem~1]{Joseph2011}.
On the other hand, if \(P=(p_i)_{i\ge1}\) is
deterministic, the infimum is not attained and
\(c_\star=1/(-\log p_{\max})\). This is the height constant of the
corresponding pebbled trie, including when there are infinitely many
possible children~\cite[Section~9]{BroutinDevroyeMcLeishDeLaSalle2008}.

\begin{theorem}[Second-order height and diameter]\label{thm:second-order}
Assume that \(\Pp(M<1)>0\) and that the infimum
in~\eqref{eq:transform-definition} is attained at \(\beta_\star\). For
\(n\ge2\) let
\[
 b_n=c_\star\log n-\frac{3c_\star}{2\beta_\star}\log\log n .
\]
Then
\begin{equation}\label{eq:second-order-tightness}
 H_n=b_n+O_{\Pp}(1),\qquad D_n=2b_n+O_{\Pp}(1).
\end{equation}
Moreover, there is a constant \(C<\infty\) such that, for all \(n\ge2\)
and \(x\ge0\),
\begin{equation}\label{eq:second-order-upper-tail}
 \Pp(D_n\ge 2b_n+2x)\le\Pp(H_n\ge b_n+x)
 \le C(1+x)\,\mathrm e^{-\beta_\star x/c_\star}.
\end{equation}
\end{theorem}

By~\eqref{eq:second-order-upper-tail},
\(\sup_{n\ge2}\E[(H_n-b_n)_+^p]<\infty\) for every \(p>0\).
The moments of the negative part depend on the lower tail of \(M\).

\begin{theorem}[Moments and lower tails]\label{thm:lower-tail}
Under the assumptions of Theorem~\ref{thm:second-order}, for every \(p>0\)
and all \(n\ge2\),
\begin{equation}\label{eq:truncated-moments}
 1+\E[(b_n-H_n)_+^p]
 \asymp 1+\E[(b_n-D_n/2)_+^p]
 \asymp 1+\E[\min\{-\log M,\log n\}^{p}],
\end{equation}
with constants depending only on \(p\) and the split law.
In particular,
\begin{equation}\label{eq:centred-moments}
 \sup_{n\ge2}\E|H_n-b_n|^p<\infty
 \quad\Longleftrightarrow\quad \E[(-\log M)^p]<\infty
 \quad\Longleftrightarrow\quad \sup_{n\ge2}\E|D_n-2b_n|^p<\infty .
\end{equation}
Each of the bounds
\begin{equation}\label{eq:lower-tail}
 \Pp(H_n\le b_n-x)\le C\,\mathrm e^{-\delta x},
 \qquad
 \Pp(D_n\le 2b_n-2x)\le C\,\mathrm e^{-\delta x}
 \qquad(n\ge2,\ x\ge0)
\end{equation}
holds for some \(C,\delta>0\) if and only if
\begin{equation}\label{eq:moment-assumption}
 \E[M^{-\eta}]<\infty\qquad\text{for some }\eta>0 .
\end{equation}
\end{theorem}

By~\eqref{eq:truncated-moments} with \(p=1\) and
\eqref{eq:second-order-upper-tail}, each of the mean relations
\(\E H_n=b_n+O(1)\) and \(\E D_n=2b_n+O(1)\) holds if and only if
\(\E[-\log M]<\infty\). If this expectation is infinite, both
\(b_n-\E H_n\) and \(2b_n-\E D_n\) are of order
\(\E[\min\{-\log M,\log n\}]\).
The variances satisfy
\begin{equation}\label{eq:variance-order}
 \Var(H_n)\asymp\Var(D_n)
 \asymp 1+\E\bigl[\min\{-\log M,\log n\}^{2}\bigr]
 \qquad(n\to\infty).
\end{equation}
If at most \(k\) coordinates are positive almost surely, then \(M\ge1/k\)
and~\eqref{eq:moment-assumption} holds.

The height is integer-valued, so \(H_n-b_n\) cannot converge in law
along the full sequence. Indeed,
\(\exp(2\pi\mathrm i(H_n-b_n))=\exp(-2\pi\mathrm i b_n)\) has more
than one limit point since \(b_n\to\infty\) and \(b_{n+1}-b_n\to0\). For binary search trees, Corre~\cite{Corre2017} showed that \(H_n-b_n\)
converges in law along every subsequence on which the fractional part
of \(b_n\) converges. We are investigating whether this also holds
under the assumptions of Theorem~\ref{thm:second-order}.
When the infimum defining \(c_\star\) is not attained, even the order
of the correction can change. For the symmetric digital search tree,
where \(P=(\tfrac12,\tfrac12,0,\dots)\) and \(c_\star\log n=\log_2n\),
the height is concentrated on two consecutive integers within bounded
distance of \(\log_2n+\sqrt{2\log_2n}-\tfrac12\log_2\log_2n\)
\cite[Theorem~3]{DrmotaFuchsHwangNeininger2021}.

\subsection{Poisson--Dirichlet splits and other examples}\label{sec:pd-heights}

For \(\theta>0\), let \(V_1,V_2,\dots\) be independent
\(\Beta(1,\theta)\) random variables and set
\begin{equation}\label{eq:gem}
 P_i=V_i\prod_{k<i}(1-V_k),\qquad i\ge1.
\end{equation}
The law of \(P\) is \(\GEM(\theta)\), and its decreasing rearrangement
has the Poisson--Dirichlet law \(\PD(0,\theta)\)~\cite{Pitman2006}.
Reordering children preserves tree shapes, so the two split laws give
the same law on sequences of unordered trees. With \(\GEM(\theta)\)
splits, labelling each occupied vertex by its ball and forgetting the
child indices gives a
recursive Ewens tree with parameter \(\theta\)
(Proposition~\ref{prop:insertion}). For both GEM laws in this section,
\(P_1\) is a size-biased pick from the ranked masses. Thus the conditional
mean of \(P_1^{\beta-1}\) given these masses is \(\sum_{i\ge1}P_i^\beta\)
for every \(\beta>0\). For \(\GEM(\theta)\), \(P_1=V_1\sim\Beta(1,\theta)\), so
\begin{equation}\label{eq:L-theta}
 L(\beta)=\E\bigl[P_1^{\beta-1}\bigr]
 =\theta\int_0^1 x^{\beta-1}(1-x)^{\theta-1}\,dx
 =\frac{\Gamma(\beta)\,\Gamma(\theta+1)}{\Gamma(\beta+\theta)}.
\end{equation}
Moreover \(p_{\max}=1\), and~\eqref{eq:moment-assumption} holds for
every \(0<\eta<1\) because
\(\max_iP_i\ge P_1\). Hence Theorems~\ref{thm:height},
\ref{thm:second-order} and~\ref{thm:lower-tail} apply. For \(\theta=1\)
we get \(L(\beta)=1/\beta\) and \(\beta_\star=c_\star=\mathrm e\), the
height constant of the uniform random recursive
tree~\cite{Devroye1987,Pittel1994}.

\begin{corollary}[Height and diameter of Plancherel trees]\label{cor:plancherel-height}
The constants \(\beta_\star>1\) and \(c_\star\) are uniquely determined by
\begin{equation}\label{eq:plancherel-critical-equation}
 \log\frac{\beta_\star(\beta_\star+1)}2=1+\frac{\beta_\star}{\beta_\star+1},
 \qquad
 c_\star=\frac{\beta_\star(\beta_\star+1)}{2\beta_\star+1}.
\end{equation}
Numerically, \(c_\star\approx 1.6738\) and
\(3c_\star/(2\beta_\star)\approx 0.8596\).
Along Fulman's upward chain \((X_n)_{n\ge1}\),
\[
 \frac{H(X_n)}{\log n}\longrightarrow c_\star,
 \qquad \frac{D(X_n)}{\log n}\longrightarrow 2c_\star
\]
almost surely and in \(L^p\) for every \(1\le p<\infty\). The same limits
hold along the insertion process of Section~\ref{sec:double-labelling}. With
\(b_n=c_\star\log n-\frac{3c_\star}{2\beta_\star}\log\log n\), the tail bounds
\eqref{eq:second-order-upper-tail} and~\eqref{eq:lower-tail} hold for
\(H(T_n^{\mathrm{Pl}})\) and \(D(T_n^{\mathrm{Pl}})\).
In particular, \(\E H(T_n^{\mathrm{Pl}})=b_n+O(1)\) and
\(\E D(T_n^{\mathrm{Pl}})=2b_n+O(1)\). Both variances are of order one.
\end{corollary}

\begin{proof}
By Proposition~\ref{prop:plancherel-growth}, both processes have the
same law on tree shapes as the split tree with \(\GEM(2)\) splits.
For this split law, \(L(\beta)=2/(\beta(\beta+1))\) by~\eqref{eq:L-theta}. The
identities~\eqref{eq:critical-constant} reduce
to~\eqref{eq:plancherel-critical-equation}. Theorems~\ref{thm:height},
\ref{thm:second-order} and~\ref{thm:lower-tail} give the stated limits,
tail bounds and expansions for the means. The variances satisfy~\eqref{eq:variance-order}.
\end{proof}

\begin{example}[Two-parameter splits]\label{ex:two-parameter}
Let \(0\le\alpha<1\) and \(\theta>-\alpha\), and let \(P\) have the
two-parameter law \(\GEM(\alpha,\theta)\), given by~\eqref{eq:gem} with
independent variables \(V_i\) of law
\(\Beta(1-\alpha,\theta+i\alpha)\)~\cite[Section~3.2]{Pitman2006}. This choice of \(P\) gives the same law on unordered rooted trees as
the Chinese restaurant tree model of Evans, Gr\"ubel and
Wakolbinger~\cite[Section~6.3]{EvansGrubelWakolbinger2012}.
The same calculation with \(P_1\sim\Beta(1-\alpha,\theta+\alpha)\) gives
\[
 L(\beta)=\frac{\Gamma(\theta+1)\,\Gamma(\beta-\alpha)}
               {\Gamma(1-\alpha)\,\Gamma(\beta+\theta)},
 \qquad\beta>\alpha.
\]
Since \(M\ge P_1\), we have \(p_{\max}=1\), and~\eqref{eq:moment-assumption}
holds for \(0<\eta<1-\alpha\). Hence
Theorems~\ref{thm:height}, \ref{thm:second-order}
and~\ref{thm:lower-tail} apply. When \(\alpha+\theta=1\), the resulting
unordered trees have the law of linear preferential attachment
trees~\cite[Theorem~1.5]{Janson2019}. In this case
\(L(\beta)=(1-\alpha)/(\beta-\alpha)\),
and~\eqref{eq:critical-constant} gives
\[
 \log\frac{\beta_\star-\alpha}{1-\alpha}
 =\frac{\beta_\star}{\beta_\star-\alpha},
 \qquad
 c_\star=\beta_\star-\alpha .
\]
For the plane-oriented recursive tree (\(\alpha=1/2\)), this gives
\(\log(2c_\star)=1+1/(2c_\star)\), in agreement with~\cite{Pittel1994}.
For the random recursive tree (\(\alpha=0\)), this gives
\(b_n=\mathrm e\log n-\tfrac32\log\log n\), as
in~\cite[Corollary~1.3]{AddarioBerryFord2013}. For \(0<\alpha<1\), \(b_n\) agrees with the centering in
\cite[Corollary~1.6]{PainSenizergues2022}, with their \((\gamma,\theta)\)
corresponding to \((1-\alpha,\beta_\star/c_\star)\).
\end{example}

\begin{example}[Binary search trees and quadtrees]\label{ex:bst}
Consider the point quadtree obtained by inserting \(n\) independent
uniform points in \([0,1]^d\), with \(d\ge1\) fixed. The split vector
consists of the volumes of the \(2^d\) boxes cut out by the coordinate
hyperplanes through the first point inserted~\cite{Devroye1987,Devroye1998}.
For \(d=1\) this is the
random binary search tree, with split vector \((U,1-U)\) for \(U\)
uniform on \([0,1]\). Independence of the coordinates gives
\(L(\beta)=(2/(\beta+1))^d\). Let \(c>2\) be the solution of
\(c\log(2\mathrm e/c)=1\). The critical equation
\eqref{eq:critical-constant} gives \(\beta_\star=c-1\) and
\(c_\star=c/d\).

Since \(M\ge2^{-d}\), Theorems~\ref{thm:second-order}
and~\ref{thm:lower-tail}, together with~\eqref{eq:variance-order}, give
\begin{align*}
 \E H_n&=\frac{c}{d}\log n-\frac{3c}{2d(c-1)}\log\log n+O(1),
 &\Var(H_n)&\asymp1,\\
 \E D_n&=2\E H_n+O(1),
 &\Var(D_n)&\asymp1.
\end{align*}
Moreover, \(\sup_{n\ge2}\E[|H_n-\E H_n|^p+|D_n-\E D_n|^p]<\infty\)
for every \(p>0\).

Devroye~\cite{Devroye1986,Devroye1987} obtained the first-order height
constant. For \(d=1\), the identity \(c/(c-1)=1/\log(c/2)\) gives
Reed's second-order term~\cite{Reed2003}. Drmota~\cite{Drmota2003}
proved that all centered absolute moments of the height are uniformly
bounded in this case.
Flajolet, Labelle, Laforest and Salvy~\cite[p.~141]{FlajoletLabelleLaforestSalvy1995}
asked about the variance of the height in higher dimensions. Majumdar, Dean and Krapivsky~\cite{MajumdarDeanKrapivsky2005} predicted, by a
non-rigorous travelling-front analysis of a continuous fragmentation model,
a \(\log\log\) correction with the same coefficient \(3c/(2d(c-1))\) and a
bounded variance of the height.
\end{example}

\subsection{Related work}\label{sec:related}

For a class of split trees with bounded branching,
Devroye~\cite[p.~412]{Devroye1998} proved an upper bound for the height
and stated the corresponding law of large numbers. Broutin, Devroye,
McLeish and de la Salle~\cite{BroutinDevroyeMcLeishDeLaSalle2008} obtained the
first-order asymptotics for a broader class of trees with a fixed
branching factor. Janson~\cite[Section~5]{Janson2019} observed that
branching random walk methods~\cite{BroutinDevroye2006,BroutinDevroyeMcLeishDeLaSalle2008}
should extend these height results to
trees with infinitely many possible children. The proof of his Theorem~5.1
gives \(\limsup_n H_n/\log n\le4/|\kappa(2)|\) almost surely.
Theorem~\ref{thm:height} identifies the limit for the split-tree model
considered here, including when the infimum defining \(c_\star\) is
not attained.

The use of branching random walks to study the height of search trees
goes back to Devroye~\cite{Devroye1986,Devroye1987}. Hu and
Shi~\cite{HuShi2009} and Addario-Berry and Reed~\cite{AddarioBerryReed2009}
identified the logarithmic correction for the minimum of a branching
random walk. A\"id\'ekon~\cite{Aidekon2013} later proved convergence
in law for the centered minimum in the non-lattice case.

For split trees, replacing the size of the subtree at \(u\) by
\(n\mu(u)\), where \(\mu(u)\) is the product of the split fractions on
the path to \(u\), gives the height approximation
\(\max\{|u|:n\mu(u)\ge1\}\). Devroye~\cite[p.~153]{Devroye2009Oberwolfach}
observed that this approximation suggests a \(\log\log n\) correction.
We use the branching random walk bounds of Madaule~\cite{Madaule2017}
and Mallein~\cite{Mallein2016}, together with the occupation-time
estimate in Lemma~\ref{lem:occupation}, to obtain the correction for
the occupied tree.

Pain and S\'enizergues~\cite{PainSenizergues2022} proved second-order
tightness for the height and diameter of preferential attachment trees
under conditions on the sequence of additive fitnesses.
Chen and Ma~\cite{ChenMa2026} study height deviations of order
\(\log n\) for uniform recursive trees.

Bertoin~\cite{Bertoin2008} and Joseph~\cite{Joseph2011} studied the
nested occupancy scheme of a multiplicative cascade, in which all
balls are sent to the leaves. Joseph found a phase transition for its
height as the capacity of the leaves varies.
For non-geometric splits with a finite minimiser \(\beta_\star\),
Joseph's results and Lemma~\ref{lem:comparison} give another proof
that \(H_n/\log n\to c_\star\) almost surely
(see Remark~\ref{rem:joseph}).

\subsection{Organization}
We prove the recursive Ewens representation of Plancherel trees
(Theorem~\ref{thm:identification}) in Section~\ref{sec:ewens-representation}.
Section~\ref{sec:nested} gives a split-tree construction of Fulman's
upward chain and establishes the branching property for general split trees.
In Section~\ref{sec:proofs-main}, we combine the first-order bounds to
prove Theorem~\ref{thm:height} on height and diameter.
For Theorem~\ref{thm:second-order}, we derive the upper tail bound in
Section~\ref{sec:upper-tail} and use occupation times to prove tightness
in Section~\ref{sec:tightness}.
Finally, Section~\ref{sec:lower-tail} proves Theorem~\ref{thm:lower-tail}
on moments and lower tails.

\subsection{Statement on AI use}
The proofs in the first version of this paper were long and indirect. Discussions with GPT-6 Astra helped the author greatly simplify the lengthy arguments. This more direct approach yielded substantially stronger results. The author also used the model to check the proofs and refine the exposition.

\subsection*{Acknowledgments}
The author thanks his PhD advisor, Pierre-Lo\"ic M\'eliot, for his
guidance and many helpful discussions, and for running simulations of
large Plancherel trees. 

\section{Plancherel trees and Ewens partitions}\label{sec:identification}

\subsection{Double increasingly labelled trees}\label{sec:double-labelling}

A \emph{double increasingly labelled tree}~\cite{KubaPanholzer2012} is a
tree equipped with an ordered pair \((\ell_1,\ell_2)\) of increasing
labellings, considered up to isomorphisms that preserve both labellings.
There are \(d(t)^2\) ordered pairs of increasing labellings on a fixed
representative of \(t\). Each orbit under \(\Aut(t)\) has size
\(|\Aut(t)|\), so the number of double increasingly labelled trees of
shape \(t\) is
\begin{equation}\label{eq:double-labelled-count}
 \frac{d(t)^2}{|\Aut(t)|}=d(t)\,u(t).
\end{equation}

These trees can be grown one vertex at a time by the insertion procedure
of Kuba and Panholzer~\cite[second proof of Theorem~4]{KubaPanholzer2012}. Start
from the one-vertex tree, whose vertex is labelled \((1,1)\). Given a
tree with \(k\) vertices, choose integers \(1\le a<b\le k+1\), attach a
new leaf \(w\) to the vertex with first label \(a\), and set
\[
 \ell_1'(v)=\ell_1(v)+\mathbf1_{\{\ell_1(v)\ge b\}},
 \quad
 \ell_2'(v)=\ell_2(v)
 \quad(v\ne w),
 \qquad
 (\ell_1'(w),\ell_2'(w))=(b,k+1).
\]
Both labellings remain increasing because \(a<b\)
(Figure~\ref{fig:double-insertion}). Conversely, the vertex \(w\)
with second label \(k+1\) is a leaf. Its first label gives \(b\),
and the first label of its parent gives \(a\).
To recover the previous tree, delete \(w\) and subtract one from each
first label greater than \(b\).
Thus the \(\prod_{k=1}^{n-1}\binom{k+1}{2}=Z_n\) sequences of choices
are in bijection with the double increasingly labelled trees with \(n\)
vertices. Independent uniform choices of the pairs give the uniform law
on these trees. In particular,
\[
 \sum_{t\in\T_n}d(t)\,u(t)=Z_n ,
\]
and the shape of a uniform double increasingly labelled tree with \(n\)
vertices has law \(p_n\).

\begin{figure}[htbp]
\centering
\begin{tikzpicture}[x=1mm,y=1mm]
\node[tree vertex] (a1) at (3.2,0) {\(1,1\)};
\node[tree vertex] (b1) at (27.6,0) {\(1,1\)};
\node[tree selected] (b2) at (27.6,10) {\(2,2\)};
\draw[tree edge] (b1)--(b2);
\node[tree vertex] (c1) at (57,0) {\(1,1\)};
\node[tree vertex] (c2) at (52,10) {\(3,2\)};
\node[tree selected] (c3) at (62,10) {\(2,3\)};
\draw[tree edge] (c1)--(c2) (c1)--(c3);
\node[tree vertex] (d1) at (91.4,0) {\(1,1\)};
\node[tree vertex] (d2) at (86.4,10) {\(3,2\)};
\node[tree vertex] (d3) at (96.4,10) {\(2,3\)};
\node[tree selected] (d4) at (86.4,20) {\(4,4\)};
\draw[tree edge] (d1)--(d2) (d1)--(d3) (d2)--(d4);
\node[tree vertex] (e1) at (125.8,0) {\(1,1\)};
\node[tree vertex] (e2) at (120.8,10) {\(4,2\)};
\node[tree vertex] (e3) at (130.8,10) {\(2,3\)};
\node[tree vertex] (e4) at (120.8,20) {\(5,4\)};
\node[tree selected] (e5) at (130.8,20) {\(3,5\)};
\draw[tree edge] (e1)--(e2) (e1)--(e3) (e2)--(e4) (e3)--(e5);
\foreach \xa/\xb/\lab in {9.4/21.4/{1,2}, 33.8/45.8/{1,2},
  68.2/80.2/{3,4}, 102.6/114.6/{2,3}}
  \draw[tree arrow] (\xa,0)--node[tree note,above=1mm]{\((\lab)\)}(\xb,0);
\foreach \x/\n in {3.2/1, 27.6/2, 57/3, 91.4/4, 125.8/5}
  \node[tree note] at (\x,-7) {\(n=\n\)};
\end{tikzpicture}
\caption{Four successive insertions, with the chosen pair \((a,b)\)
above each arrow. The new vertex is shaded.}
\label{fig:double-insertion}
\end{figure}

\subsection{The Ewens representation}\label{sec:ewens-representation}

Fix an increasingly labelled tree of shape \(t\in\T_n\). At each vertex
\(v\), the label sets of the child subtrees form a partition of the
\(|t_v|-1\) labels below \(v\). Multiplying the Ewens probabilities~\eqref{eq:ewens} of
these partitions gives the probability \(w_\theta(t)\) of this labelled tree,
\begin{equation}\label{eq:general-labelled}
 w_\theta(t)
 =\frac{\theta^{\,n-1}\prod_{v\ne\mathrm{root}}(|t_v|-1)!}
        {\prod_{v\in t}\theta^{\overline{|t_v|-1}}}.
\end{equation}
For \(\theta=1\), all factors cancel except the denominator at the root.
Thus \(w_1(t)=1/(n-1)!\) for every increasingly labelled tree, which
gives the uniform random recursive tree.

\begin{proof}[Proof of Theorem~\ref{thm:identification}]
Since \(2^{\overline m}=(m+1)!\),
\begin{equation}\label{eq:labelled-probability}
 w_2(t)
 =\frac{2^{n-1}\prod_{v\ne\mathrm{root}}(|t_v|-1)!}
       {\prod_{v\in t}|t_v|!}
 =\frac{2^{n-1}}{(n-1)!\prod_{v\in t}|t_v|}
 =\frac{d(t)}{Z_n}
\end{equation}
by~\eqref{eq:hook} and~\eqref{eq:plancherel}. Summing over the \(u(t)\)
increasingly labelled trees of shape \(t\) yields
\(\Pp(T_n^{(2)}=t)=p_n(t)\).
\end{proof}

\subsection{Comparison with Young diagrams}\label{sec:comparison}

Table~\ref{tab:comparison} compares the Plancherel measures on
Young diagrams and on rooted trees. Here \(h(\square)\) is the hook length of a cell, and \(U\)
and \(D\) are the up and down operators. On rooted trees, \(U\) attaches a leaf at each vertex and \(D\) removes
each non-root leaf. Different choices are counted separately, even when
they yield isomorphic trees. Thus \(u(t)\) counts
the upward paths from the one-vertex tree to \(t\), and \(d(t)\) counts
the downward paths back to the one-vertex tree. The commutation relations are due
to Stanley~\cite{Stanley1988} and
Hoffman~\cite[Proposition~2.2]{Hoffman2003}. Like a recording tableau, the second labelling in
Section~\ref{sec:double-labelling} gives the insertion time of each vertex. The resulting growth chain is
identified in Proposition~\ref{prop:plancherel-growth}. For a vertex
\(v\) of \(t\), we write \(t^{+v}\) for the tree obtained by attaching a
new leaf to \(v\).

\begin{table}[htbp]
\centering
\footnotesize
\begin{tabular}{@{}lll@{}}
\toprule
 & partitions \(\lambda\) of \(n\) & rooted trees \(t\in\T_n\)\\
\midrule
hook-length formula
 & \(f^\lambda=n!/\prod_{\square\in\lambda}h(\square)\)
 & \(d(t)=n!/\prod_{v\in t}|t_v|\)\\
upward and downward paths
 & \(f^\lambda\) and \(f^\lambda\)
 & \(u(t)\) and \(d(t)\)\\
commutation relation
 & \(DU-UD=I\)
 & \(DU-UD=nI\)\\
Plancherel measure
 & \((f^\lambda)^2/n!\)
 & \(d(t)u(t)/Z_n\)\\
bijection
 & Robinson--Schensted
 & Kuba--Panholzer\\
probability per upward edge
 & \(f^\nu/((n+1)f^\lambda)\)
 & \(d(t')/(\binom{n+1}{2}d(t))\)\\
\bottomrule
\end{tabular}\par\medskip
\caption{Plancherel measures on partitions and on rooted trees. The commutation relation
for rooted trees holds on the span of trees with \(n\) vertices.
In the last row, \(\nu\) is obtained from \(\lambda\) by adding a box
and \(t'=t^{+v}\). The last entry gives the probability of attaching a leaf
at a fixed vertex \(v\). Summing over all \(v\) with \(t^{+v}\cong t'\)
gives the transition probability from \(t\) to \(t'\).}
\label{tab:comparison}
\end{table}

\section{Split trees in a random environment}\label{sec:nested}

\subsection{The nested construction}\label{sec:coupling}

Let \(\U=\bigcup_{h\ge0}\{1,2,\dots\}^h\) be the Ulam tree, with root
\(\varnothing\). For \(u\in\U\) we write \(|u|\) for its depth, \(ui\)
for its \(i\)th child, \(u|_k\) for its ancestor at depth \(k\le|u|\),
and \(u\,\U=\{uv:v\in\U\}\) for the subtree rooted at \(u\). We write
\(u\preceq v\) if \(v\in u\,\U\). Let \((P(u))_{u\in\U}\) be independent
copies of the split vector \(P\), let \(\cE\) be the \(\sigma\)-field
that they generate, and define the masses
\[
 \mu(\varnothing)=1,
 \qquad
 \mu(ui)=\mu(u)\,P_i(u),
\]
so that \(\sum_{|u|=h}\mu(u)=1\) for every \(h\ge0\). Conditionally on \(\cE\), assign independent random rays \((R_b)_{b\ge1}\)
to the balls in order of insertion. Each ray starts at the root and
moves from \(u\) to \(ui\) with probability \(P_i(u)\). Thus
\begin{equation}\label{eq:ray-law}
 \Pp(u\in R_b\mid\cE)=\mu(u).
\end{equation}
Ball \(b\) follows \(R_b\) and is stored at the first unoccupied vertex
\(s_b\) of \(R_b\). The set
\(\widehat T_n=\{s_1,\dots,s_n\}\) is a tree with \(n\) vertices, in which
the ball labels increase along every path from the root. Moreover
\(\widehat T_n\subset\widehat T_{n+1}\), so that \(H_n=H(\widehat T_n)\)
is non-decreasing. This realisation of the split tree in
Section~\ref{sec:height-intro} follows~\cite[proof of Theorem~5.1]{Janson2019}.

\subsection{Growth of recursive Ewens trees}\label{sec:pd-splits}

\begin{proposition}\label{prop:insertion}
Let \(P\) have the \(\GEM(\theta)\) law~\eqref{eq:gem}. For each \(n\ge1\), keep the ball labels on \(\widehat T_n\) and forget
the child indices. The resulting tree has the law of a recursive Ewens
tree with parameter \(\theta\) on \(\{1,\dots,n\}\).
\end{proposition}

\begin{proof}
We argue by induction on the number \(n\) of balls, allowing any set of
\(n\) labels inserted in increasing order. The smallest label is stored at
the root. Conditionally on \(P(\varnothing)\), the other balls choose a
child of the root independently with law \(P(\varnothing)\). After
averaging over \(P(\varnothing)\), the partition of these labels by their
chosen child has law \(\Ewens(n-1,\theta)\) by Kingman's paintbox
representation~\cite{Pitman2006}. Given this partition, the balls in
each block enter the corresponding subtree in increasing order. Each
subtree has an independent copy of the environment. By induction, the
subtrees are independent recursive Ewens trees on their respective label sets.
\end{proof}

In particular, every increasingly labelled tree of shape \(t\) has
probability \(w_\theta(t)\), given by~\eqref{eq:general-labelled}.
For a vertex \(v\) of \(t\in\T_n\), let
\(\rho_0,\dots,\rho_r=v\) be the path from the root \(\rho_0\) to \(v\),
and set
\[
 q_\theta(t,v)=\frac{\theta}{\theta+|t_v|-1}
 \prod_{i=0}^{r-1}\frac{|t_{\rho_{i+1}}|}{\theta+|t_{\rho_i}|-1}.
\]
This is the probability that a new leaf is attached to \(v\) under the
Chinese restaurant rule. At an occupied vertex \(w\), the new ball
creates a child of \(w\) with probability \(\theta/(\theta+|t_w|-1)\).
It enters the subtree of a given existing child \(w'\) with probability
\(|t_{w'}|/(\theta+|t_w|-1)\). For \(\theta=1\) the product telescopes
and \(q_1(t,v)=1/|t|\).

Fulman's upward chain~\cite[Definition~2]{Fulman2009} starts from the
one-vertex tree and has transition probabilities
\begin{equation}\label{eq:fulman-chain}
 \Pp(X_{n+1}=t'\mid X_n=t)
 =\frac{d(t')}{\binom{n+1}{2}\,d(t)}\,\#\{v\in t:t^{+v}\cong t'\},
\end{equation}
for \(t\in\T_n\) and \(t'\in\T_{n+1}\), where \(v\) runs over the
vertices of a representative of \(t\) and \(\cong\) denotes isomorphism of rooted trees. 

\begin{proposition}\label{prop:plancherel-growth}
Let \(P\) have the \(\GEM(\theta)\) law. The shapes of
\((\widehat T_n)_{n\ge1}\) form a Markov chain started from the one-vertex
tree, with transition probabilities
\begin{equation}\label{eq:plancherel-growth}
 \Pp\bigl(\widehat T_{n+1}\cong t'\mid\widehat T_n\cong t\bigr)
 =\sum_{v\in t:\,t^{+v}\cong t'}q_\theta(t,v).
\end{equation}
For \(\theta=2\) this is Fulman's upward chain~\eqref{eq:fulman-chain},
and it is also the chain of shapes in the insertion process of
Section~\ref{sec:double-labelling}.
\end{proposition}

\begin{proof}
\emph{Step 1. Transitions between tree shapes.}
In this step we view \(\widehat T_n\) as an unordered tree with its arrival labels. Let \(\mathbf t\) be an increasingly labelled
tree of shape \(t\in\T_n\), and let \(\mathbf t^{+v}\) be obtained by
attaching a leaf labelled \(n+1\) to \(v\). Since the labelled tree
\(\widehat T_{n+1}\) determines \(\widehat T_n\) by restriction,
Proposition~\ref{prop:insertion} and~\eqref{eq:general-labelled} give
\[
 \Pp(\widehat T_{n+1}=\mathbf t^{+v}\mid\widehat T_n=\mathbf t)
 =\frac{w_\theta(t^{+v})}{w_\theta(t)}
 =\frac{\theta\prod_{i=1}^r|t_{\rho_i}|}
        {\prod_{i=0}^r(\theta+|t_{\rho_i}|-1)}
 =q_\theta(t,v),
\]
since the subtree sizes change only along the path
\(\rho_0,\dots,\rho_r\). The labelled tree \(\mathbf t\) determines all earlier shapes by restriction.
For fixed \(t'\), the sum of \(q_\theta(t,v)\) over vertices \(v\) with
\(t^{+v}\cong t'\) depends on \(\mathbf t\) only through its shape \(t\). Hence the shapes form a Markov chain with the
transition probabilities~\eqref{eq:plancherel-growth}.

\smallskip\noindent\emph{Step 2. The two growth processes at \(\theta=2\).}
By~\eqref{eq:labelled-probability} and \(Z_{n+1}/Z_n=\binom{n+1}{2}\),
\[
 q_2(t,v)=\frac{w_2(t^{+v})}{w_2(t)}
 =\frac{d(t^{+v})}{\binom{n+1}{2}d(t)}.
\]
Summing over \(v\) with \(t^{+v}\cong t'\) gives~\eqref{eq:fulman-chain}.
In the insertion process of Section~\ref{sec:double-labelling},
restricting the second labels to \(1,\dots,k\) recovers the tree at
time \(k\) with its second labelling. At time \(n\), each
\(\ell_2\)-labelled tree of shape \(t\) admits \(d(t)\) choices for
\(\ell_1\) and hence has probability \(d(t)/Z_n=w_2(t)\). The probability
of attaching at \(v\) is therefore \(q_2(t,v)\), and summing over \(v\)
gives~\eqref{eq:plancherel-growth} for the tree shapes.
\end{proof}

\subsection{Occupancy counts and the branching property}\label{sec:occupancy}

We return to a general split vector. For \(u\in\U\) and \(n\ge1\) let
\[
 N_n(u)=\#\{b\le n:u\in R_b\},
 \qquad
 K_n(u)=\#\{b\le n:u\preceq s_b\}.
\]
Among the first \(n\) balls, \(N_n(u)\) counts those whose rays pass
through \(u\), even if the ball stops at a strict ancestor. The count
\(K_n(u)\) includes only those stored at \(u\) or one of its descendants.
In particular, \(u\in\widehat T_n\) if and only if \(K_n(u)\ge1\). By~\eqref{eq:ray-law}, conditionally on \(\cE\),
\begin{equation}\label{eq:binomial-counts}
 N_n(u)\sim\Bin(n,\mu(u)).
\end{equation}

\begin{lemma}[Comparison]\label{lem:comparison}
For \(n\ge1\) and \(u\in\U\) we have \(0\le N_n(u)-K_n(u)\le|u|\).
Almost surely, \(N_n(u)/n\to\mu(u)\) and \(K_n(u)/n\to\mu(u)\)
simultaneously for all \(u\in\U\).
If \(H_n\ge h+j-1\) for integers \(h\ge0\) and \(j\ge1\), then
\(N_n(u)\ge j\) for some \(u\) at depth \(h\).
\end{lemma}

\begin{proof}
Since \(s_b\in R_b\),
\[
 N_n(u)-K_n(u)
 =\#\{b\le n:u\in R_b,\ |s_b|<|u|\}\le|u|,
\]
because each strict ancestor of \(u\) stores at most one ball
(Figure~\ref{fig:occupancy}). Given \(\cE\), the indicators
\(\mathbf1_{\{u\in R_b\}}\), \(b\ge1\), are i.i.d. with mean \(\mu(u)\).
The strong law gives \(N_n(u)/n\to\mu(u)\) simultaneously for all
\(u\in\U\), since \(\U\) is countable. The comparison gives the same
limit for \(K_n(u)/n\).
For the last assertion, take a vertex of \(\widehat T_n\) at depth
\(h+j-1\) and its ancestor \(u\) at depth \(h\).
The \(j\) vertices of the path between them belong to \(\widehat T_n\),
so that \(N_n(u)\ge K_n(u)\ge j\).
\end{proof}

\begin{figure}[htbp]
\centering
\begin{tikzpicture}[x=1mm,y=1mm]
\node[tree vertex] (o) at (0,0) {1};
\node[tree vertex] (a) at (-5,10) {2};
\node[tree vertex] (b) at (5,10) {3};
\node[tree selected] (u) at (-10,20) {4};
\node[tree vertex] (c) at (0,20) {5};
\node[tree selected] (d) at (-15,30) {6};
\node[tree selected] (e) at (-5,30) {7};
\draw[tree edge] (o)--(a) (o)--(b) (a)--(u) (a)--(c) (u)--(d) (u)--(e);
\draw[tree ray] (d)--(-20,40) node[tree note,above]{\(R_1\)};
\draw[tree ray] (e)--(0,40) node[tree note,above]{\(R_2\)};
\node[tree note,below=1mm] at (o.south) {\(\varnothing\)};
\node[tree note,left=1mm] at (u.west) {\(u\)};
\node[tree note,anchor=west,align=left] at (20,17) {%
  \(\{b\le7:u\in R_b\}=\{1,2,4,6,7\}\)\\[3pt]
  \(\{b\le7:u\preceq s_b\}=\{4,6,7\}\)};
\end{tikzpicture}
\caption{The tree \(\widehat T_7\) with each vertex labelled by the ball
that it stores, and the subtree of \(u=s_4\) shaded. The rays \(R_1\) and
\(R_2\) pass through \(u\), but balls \(1\) and \(2\) are stored at its two
strict ancestors, so that \(N_7(u)-K_7(u)=5-3=|u|\).}
\label{fig:occupancy}
\end{figure}

By Lemma~\ref{lem:comparison}, almost surely \(K_n(w)\to\infty\) for
every vertex \(w\) of positive mass. Vertices of zero mass receive no
balls, by~\eqref{eq:binomial-counts} and \(K_n(w)\le N_n(w)\).
For \(\mu(w)>0\) and \(m\ge1\), let \(n\) be the first time at which
\(K_n(w)=m\), and set
\[
 \widehat T^{\,w}_m=\{u\in\U:wu\in\widehat T_n\}.
\]
This is the tree formed by the first \(m\) balls stored in \(w\,\U\),
under the identification \(wu\mapsto u\).

Let \(w_1,w_2\) be incomparable vertices, meaning that neither is an
ancestor of the other. If \(1\le m_i\le K_n(w_i)\) for \(i=1,2\),
the path joining deepest vertices of the two subtrees passes through
\(w_1\) and \(w_2\). Hence
\begin{equation}\label{eq:two-subtree-diameter}
 H(\widehat T^{\,w_1}_{m_1})+\operatorname{dist}(w_1,w_2)
 +H(\widehat T^{\,w_2}_{m_2})\le D_n\le2H_n ,
\end{equation}
where \(\operatorname{dist}\) denotes graph distance in \(\U\).

\begin{lemma}[Branching property]\label{lem:branching}
Let \(w_1,\dots,w_r\) be pairwise incomparable random vertices. Assume that
\(\{w_i=w\}\in\sigma(P(w|_k):0\le k<|w|)\) for every \(i\) and
\(w\in\U\). Let \(\mathcal G\) be the
\(\sigma\)-field generated by \(w_1,\dots,w_r\), by the split vectors at vertices outside \(\bigcup_iw_i\,\U\), and by the counts
\(K_n(w_i)\) for \(n\ge1\) and \(1\le i\le r\). Conditionally on
\(\mathcal G\), the subtree processes \((\widehat T^{\,w_i}_m)_{m\ge1}\)
at the vertices \(w_i\) of positive mass are independent copies of
\((\widehat T_m)_{m\ge1}\).
\end{lemma}

\begin{proof}
\emph{Step 1. Balls entering each subtree.}
Assume first that the \(w_i\) are deterministic, and let \(R'_b\) be the
ray \(R_b\) stopped at the first of the vertices \(w_1,\dots,w_r\) that
it visits, if any. Let \(B_i=\{b:s_b\in w_i\,\U\}\). By the insertion rule,
\[
 b\in B_i\quad\Longleftrightarrow\quad
 w_i\in R'_b\ \text{ and }\
 \{w_i|_k:0\le k<|w_i|\}\subseteq\{s_1,\dots,s_{b-1}\}.
\]
Following the stopped rays in order determines the occupied vertices
outside the chosen subtrees. Each ball is stored at the first empty
vertex it reaches there or is assigned to \(B_i\) upon reaching \(w_i\).
This determines all \(B_i\) and \(K_n(w_i)=\#(B_i\cap\{1,\dots,n\})\).

\smallskip\noindent\emph{Step 2. Conditional law of the subtrees.}
Given \(\cE\) and the stopped rays \((R_b')_{b\ge1}\), the continuations below \(w_i\) for
\(b\in B_i\), \(1\le i\le r\), are independent, with transition
probabilities given by \((P(w_iu))_{u\in\U}\). In arrival order they
generate the subtree process at \(w_i\) by the insertion rule.
The stopped rays use only splits outside
\(\bigcup_iw_i\,\U\). Given those splits and the stopped rays, the
environments inside the \(w_i\,\U\) are independent copies of the
original environment. The outside splits determine which \(w_i\) have
positive mass. For these vertices, \(B_i\) is infinite and the subtree
processes have the stated joint law given the outside splits and the
stopped rays, and hence also given \(\mathcal G\).

Now suppose that the \(w_i\) are random. For fixed pairwise incomparable
vertices \(w'_1,\dots,w'_r\), the event \(\{w_1=w'_1,\dots,w_r=w'_r\}\)
is determined by the splits at the strict ancestors of the \(w'_i\).
These ancestors lie outside \(\bigcup_i w'_i\,\U\), so the deterministic case applies on this event.
\end{proof}

\section{First-order height and diameter}\label{sec:height}

\subsection{The functions \texorpdfstring{\(L\)}{L} and
\texorpdfstring{\(\kappa\)}{kappa}}\label{sec:transform}

For each integer \(h\ge0\), let
\[
 \mathcal F_h=\sigma\bigl(P(v):|v|<h\bigr),
 \qquad
 W_h(\beta)=\sum_{|u|=h}\mu(u)^\beta,\quad \beta>0 .
\]
Recall that \(L(\beta)=\E[\sum_{i\ge1}P_i^\beta]=\E W_1(\beta)\)
and \(\kappa(\beta)=\log L(\beta)\).

\begin{lemma}\label{lem:transform}
Part~\ref{it:moments} holds for every split vector. For the other parts,
assume that \(\Pp(M<1)>0\).
\begin{enumerate}[label=\textup{(\roman*)},leftmargin=*]
\item\label{it:L} \(L(1)=1\), and \(L\) is non-increasing on \((0,\infty)\)
and strictly decreasing on \([1,\infty)\). For \(\beta>1\),
\(M^\beta\le\sum_iP_i^\beta\le M^{\beta-1}\) and \(0<L(\beta)<1\).
\item\label{it:kappa} \(\kappa\) is convex on \([1,\infty)\) and smooth
on \((1,\infty)\), and \(-1/\kappa(2)\le c_\star<\infty\).
\item\label{it:moments} For \(\beta>0\), \(\E W_h(\beta)=L(\beta)^h\)
for every \(h\ge0\). If \(L(\beta)<\infty\), the normalised sums
\(L(\beta)^{-h}W_h(\beta)\) form a non-negative
\((\mathcal F_h)\)-martingale with mean one.
\item\label{it:attained} \(L(\beta)^{1/\beta}\to p_{\max}\) as
\(\beta\to\infty\), and~\eqref{eq:attainment} holds.
If the infimum in~\eqref{eq:transform-definition} is attained, then
\(1/c_\star>-\log p_{\max}\). The minimiser \(\beta_\star\in(1,\infty)\)
is unique and satisfies
\begin{equation}\label{eq:critical-constant}
 \kappa(\beta_\star)=\beta_\star\,\kappa'(\beta_\star)
 =-\frac{\beta_\star}{c_\star}.
\end{equation}
\end{enumerate}
\end{lemma}

\begin{proof}
\ref{it:L} The sum \(\sum_iP_i^\beta\) is non-increasing in \(\beta\)
and equals one at \(\beta=1\). For \(1\le a<b\), we have
\(\sum_iP_i^b<\sum_iP_i^a\le1\) on \(\{M<1\}\), so \(L(b)<L(a)\).
For \(\beta>1\), the largest term is \(M^\beta>0\) and
\(\sum_iP_i^\beta\le M^{\beta-1}\sum_iP_i=M^{\beta-1}\).

\ref{it:kappa} Convexity follows from H\"older's inequality. For every
integer \(r\ge0\), the function \(x\mapsto x^{\beta-1}|\log x|^r\) is
bounded on \((0,1]\), locally uniformly in \(\beta>1\). Since
\(\sum_iP_i=1\), differentiation under the expectation shows that
\(L\), and hence \(\kappa\), is smooth on \((1,\infty)\).
By~\ref{it:L} and Jensen's inequality, for every \(\beta>1\),
\(L(\beta)\ge\E[M^\beta]\ge(\E M)^\beta\ge L(2)^\beta\).
This gives \(c_\star\ge-1/\kappa(2)\), and taking \(\beta=2\)
in~\eqref{eq:transform-definition} gives \(c_\star\le-2/\kappa(2)<\infty\).

\ref{it:moments} For \(|u|=h\), the mass \(\mu(u)\) is
\(\mathcal F_h\)-measurable and \(P(u)\) is independent of \(\mathcal F_h\).
Summing over the children of each vertex and using Tonelli's theorem gives
\[
 \E[W_{h+1}(\beta)\mid\mathcal F_h]
 =\sum_{\substack{|u|=h\\\mu(u)>0}}\mu(u)^\beta\,
    \E\!\left[\sum_i P_i(u)^\beta\right]
 =L(\beta)W_h(\beta).
\]
Taking expectations and iterating from \(W_0(\beta)=1\) gives the mean
formula. When \(L(\beta)<\infty\), dividing by \(L(\beta)^{h+1}\)
gives the martingale property.

\ref{it:attained} Write \(p=p_{\max}\), \(\nu=\E\,\#\{i:P_i=p\}\) and
\(f(\beta)=-\kappa(\beta)/\beta\), so that \(1/c_\star=\sup_{\beta>1}f(\beta)\).
Since \(0\le P_i/p\le1\) and \(\sum_iP_i/p=1/p\), dominated convergence gives
\begin{equation}\label{eq:normalised-L}
 L(\beta)\,p^{-\beta}
 =\E\sum_{i\ge1}\Bigl(\frac{P_i}{p}\Bigr)^{\beta}
 \downarrow\nu\qquad(\beta\to\infty).
\end{equation}
The bounds \(\E M^\beta\le L(\beta)\le p^{\beta-1}\) from part~\ref{it:L}
and \(\|M\|_\beta\to\|M\|_\infty=p\) give
\(L(\beta)^{1/\beta}\to p\). Thus \(f(\beta)\to-\log p\) as
\(\beta\to\infty\), and dominated convergence at \(\beta=1\) gives
\(L(\beta)\to1\) and hence \(f(\beta)\to0\) as \(\beta\downarrow1\).

If \(\nu<1\), then \eqref{eq:normalised-L} gives a \(\beta_0>1\) with
\(f(\beta_0)>-\log p\ge0\). The continuous function \(f\) therefore
attains its supremum in \((1,\infty)\). Now suppose that \(\nu\ge1\). Then \(L(\beta)p^{-\beta}\ge\nu\ge1\),
so \(f(\beta)\le-\log p\). Since \(f(\beta)\to-\log p\) as
\(\beta\to\infty\), we have \(\sup_{\beta>1}f(\beta)=-\log p\). If the infimum in~\eqref{eq:transform-definition} were attained at some
\(\beta>1\), then \(f(\beta)=-\log p\), hence \(L(\beta)p^{-\beta}=1\).
This would imply \(\nu=1\) and
\(\E\sum_{i:\,0<P_i<p}(P_i/p)^\beta=0\).
All positive coordinates would then equal \(p\), giving
\(\nu=\E\,\#\{i:P_i=p\}=1/p\). Since \(\nu=1\), this would imply
\(M=p=1\) almost surely, contradicting \(\Pp(M<1)>0\). This proves~\eqref{eq:attainment}.

Suppose now that the infimum defining \(c_\star\) is attained, so that
\(\nu<1\). Differentiating \(\kappa=\log L\) gives, for \(\beta>1\),
\begin{equation}\label{eq:kappa-variance}
 \kappa''(\beta)=\frac{1}{L(\beta)}\,
   \E\sum_{i:P_i>0}P_i^\beta\bigl(\log P_i-\kappa'(\beta)\bigr)^2>0 .
\end{equation}
Indeed, if this variance vanished, all positive coordinates would equal
the same constant, necessarily \(p\), giving \(\nu=1/p\ge1\).
Thus \(\beta^2 f'(\beta)=\kappa(\beta)-\beta\kappa'(\beta)\) has
derivative \(-\beta\kappa''(\beta)<0\). It follows that \(f'\) has at
most one zero, so the maximiser \(\beta_\star\) is unique. The identities
\(f'(\beta_\star)=0\) and \(f(\beta_\star)=1/c_\star\)
give~\eqref{eq:critical-constant}.
\end{proof}

Given \(\mathcal F_k\), the families of relative masses
\((\mu(vw)/\mu(v))_{w\in\U}\), indexed by the vertices \(v\) of positive
mass at depth \(k\), are independent copies of \((\mu(w))_{w\in\U}\).
Hence, for each such \(v\), every integer \(r\ge0\) and all \(\beta,q>0\),
\begin{equation}\label{eq:subtree-moments}
 \E\!\left[\left(\sum_{|w|=r}\mu(vw)^\beta\right)^q
     \,\middle|\,\mathcal F_k\right]
 =\mu(v)^{\beta q}\,\E[W_r(\beta)^q].
\end{equation}
For \(q=1\), the right-hand side is \(\mu(v)^\beta L(\beta)^r\)
by Lemma~\ref{lem:transform}\ref{it:moments}.

\subsection{Upper bound}\label{sec:upper}

\begin{proposition}\label{prop:upper}
For integers \(n\ge1\), \(h\ge0\), \(j\ge1\) and every \(0<\beta\le j\),
\begin{equation}\label{eq:quenched-tail}
 \Pp(H_n\ge h+j-1\mid\cE)
 \le\min\!\left\{1,\frac{n^\beta\,W_h(\beta)}{(j!)^{\beta/j}}\right\},
 \qquad
 \Pp(H_n\ge h+j-1)\le\frac{n^\beta L(\beta)^h}{(j!)^{\beta/j}}.
\end{equation}
If moreover \(\Pp(M<1)>0\), then for \(\beta>1\),
\(j=\lceil\beta\rceil\) and
\(a_\beta=\beta/(-\kappa(\beta))\),
\begin{equation}\label{eq:exponential-upper-tail}
 \Pp(H_n\ge a_\beta\log n+x+j-1)\le L(\beta)^x,
 \qquad n\ge1,\ x\ge0,
\end{equation}
and \(\limsup_{n\to\infty}H_n/\log n\le c_\star\) almost surely.
\end{proposition}

\begin{proof}
By Lemma~\ref{lem:comparison}, a union bound
and~\eqref{eq:binomial-counts},
\[
 \Pp(H_n\ge h+j-1\mid\cE)
 \le\sum_{|u|=h}\min\Bigl\{1,\binom nj\mu(u)^j\Bigr\}
 \le\sum_{|u|=h}\min\Bigl\{1,\frac{(n\mu(u))^j}{j!}\Bigr\}.
\]
Using \(\min\{1,y\}\le y^{\beta/j}\) in the sum and bounding the
conditional probability by \(1\) gives the first inequality
in~\eqref{eq:quenched-tail}. The second follows by taking expectations
and applying Lemma~\ref{lem:transform}\ref{it:moments}.

Assume now that
\(\Pp(M<1)>0\), and let \(\beta>1\) and
\(h=\lceil a_\beta\log n+x\rceil\). As \(H_n\) is an integer, the event
in~\eqref{eq:exponential-upper-tail} is \(\{H_n\ge h+j-1\}\), and
\(n^\beta L(\beta)^h\le n^\beta L(\beta)^{a_\beta\log n+x}=L(\beta)^x\).
Fix \(\varepsilon>0\). By~\eqref{eq:exponential-upper-tail} with
\(n=2^k\) and \(x=\varepsilon k\), and by the Borel--Cantelli lemma,
almost surely \(H_{2^k}\le(a_\beta\log2+\varepsilon)k+j\) for all large
\(k\). For \(2^{k-1}\le n<2^k\) we have \(H_n\le H_{2^k}\) and
\(\log n\ge(k-1)\log2\), so that
\(\limsup_nH_n/\log n\le a_\beta+\varepsilon/\log2\). Let \(\varepsilon\downarrow0\) and take the infimum over rational
\(\beta>1\). By Lemma~\ref{lem:transform}\ref{it:kappa},
\(\beta\mapsto a_\beta\) is continuous on \((1,\infty)\), so this
infimum equals \(c_\star\).
\end{proof}

\subsection{Lower bound}\label{sec:lower}

For \(h\ge0\) let \(\mu_h^\star=\max_{|u|=h}\mu(u)\).

\begin{lemma}\label{lem:largest-mass}
Assume that \(\Pp(M<1)>0\). Then, almost surely,
\begin{equation}\label{eq:largest-mass-speed}
 \lim_{h\to\infty}\frac{-\log\mu_h^\star}{h}=\frac1{c_\star}.
\end{equation}
\end{lemma}

\begin{proof}
Place a particle at \(\log\mu(u)\) for each vertex \(u\) with \(\mu(u)>0\).
Since \(\mu(ui)=\mu(u)P_i(u)\), these particles form a branching random
walk with reproduction point process \(\Xi=\sum_{i:P_i>0}\delta_{\log P_i}\). Every particle has at least one
child and, with positive probability, at least two. The log-Laplace
transform \(\log\E\int\mathrm e^{\beta x}\,\Xi(dx)=\kappa(\beta)\)
satisfies \(\kappa(1)=\log\E\sum_iP_i=0\). Biggins' theorem on the speed of the rightmost
particle allows infinitely many children~\cite[Section~2 and Theorem~3.1]{Biggins2010}
and gives
\(h^{-1}\log\mu_h^\star\to\inf_{\beta>0}\kappa(\beta)/\beta\) almost
surely. Since \(\kappa\ge0\) on \((0,1]\), we have
\(\inf_{\beta>0}\kappa(\beta)/\beta
=\inf_{\beta>1}\kappa(\beta)/\beta=-1/c_\star\)
by~\eqref{eq:transform-definition}.
\end{proof}

\begin{proposition}\label{prop:lower}
Assume that \(\Pp(M<1)>0\). Then
\(\liminf_{n\to\infty}H_n/\log n\ge c_\star\) almost surely.
\end{proposition}
\begin{proof}
Fix \(h\ge0\). If \(H_n<h\), then \(K_n(u)=0\) for every vertex \(u\)
at depth \(h\), so \(N_n(u)\le h\) by Lemma~\ref{lem:comparison}.
Given \(\cE\), \(N_n(u)\) has law \(\Bin(n,\mu(u))\)
by~\eqref{eq:binomial-counts}. Applying the Chernoff bound
at a vertex of maximal mass gives
\begin{equation}\label{eq:lower-chernoff}
 \Pp(H_n<h\mid\cE)
 \le \exp\Bigl(-\frac{n\mu_h^\star}{8}\Bigr)
 \qquad\text{if }n\mu_h^\star\ge2h.
\end{equation}
Fix \(c\in(0,c_\star)\) and set \(h_n=\lfloor c\log n\rfloor\).
By~\eqref{eq:largest-mass-speed},
\(n\mu_{h_n}^\star=n^{1-c/c_\star+o(1)}\) almost surely.
Since \(1-c/c_\star>0\) and \(h_n=O(\log n)\),
\eqref{eq:lower-chernoff} gives
\(\sum_{n\ge2}\Pp(H_n<h_n\mid\cE)<\infty\) almost surely.
By the conditional Borel--Cantelli lemma, \(H_n\ge h_n\) eventually
almost surely. This holds simultaneously for all rational
\(c\in(0,c_\star)\), so \(\liminf_{n\to\infty}H_n/\log n\ge c_\star\)
almost surely.
\end{proof}

\subsection{Height and diameter}\label{sec:proofs-main}

\begin{proof}[Proof of Theorem~\ref{thm:height}]
\emph{Step 1. Height and moment bounds.}
Propositions~\ref{prop:upper} and~\ref{prop:lower} give the almost sure
height limit. By Lemma~\ref{lem:transform}\ref{it:moments} and the
martingale convergence theorem, \(W=\sup_{h\ge0}L(2)^{-h}W_h(2)\)
is finite almost surely. For
\(n\ge2\), \(x\ge0\) and \(h=\lceil a_2\log n+x\rceil\),
\eqref{eq:quenched-tail} with \(\beta=j=2\) gives
\[
 \Pp(H_n\ge a_2\log n+x+1\mid\cE)
 \le\tfrac12\,W\,L(2)^x,
 \qquad
 \Pp(H_n\ge a_2\log n+x+1)\le\tfrac12\,L(2)^x .
\]
For each integer \(q\ge2\), integrating the conditional tail bound
against \(q x^{q-1}\,dx\) gives
\(\E[(H_n-a_2\log n-1)_+^q\mid\cE]
 \le \tfrac{qW}{2}\int_0^\infty x^{q-1}L(2)^x\,dx\).
The unconditional bound is the same with \(W\) replaced by \(1\).
Since \(L(2)<1\), we obtain
\[
 \sup_{n\ge2}\E[(H_n/\log n)^q]<\infty,
 \qquad
 \sup_{n\ge2}\E[(H_n/\log n)^q\mid\cE]<\infty
 \quad\text{almost surely}.
\]
For almost every environment, the conditional bounds hold for all
integers \(q\ge2\), and \(H_n/\log n\to c_\star\) almost surely under
the conditional law. For any \(p\ge1\), choosing an integer \(q>p\)
gives uniform integrability of \(((H_n/\log n)^p)_{n\ge2}\).
The almost sure convergence therefore yields \(L^p\) convergence both
unconditionally and for each such environment.

\smallskip\noindent\emph{Step 2. Diameter.}
The tree of positive-mass vertices branches at finite depth almost surely,
since it has only one vertex at depth \(k\) with probability
\(\Pp(M=1)^k\to0\). At its first branching vertex, let \(w,w'\) be
the two positive-mass children with the smallest indices.
For \(v\in\{w,w'\}\), \(K_n(v)/n\to\mu(v)>0\) almost surely by Lemma~\ref{lem:comparison}.
By Lemma~\ref{lem:branching} and the height limit just proved,
\(H(\widehat T^{\,v}_m)/\log m\to c_\star\) almost surely as the number
\(m\) of balls stored in the subtree rooted at \(v\) tends to infinity.
For all large \(n\), apply~\eqref{eq:two-subtree-diameter} with
\(m_1=K_n(w)\) and \(m_2=K_n(w')\). Since \(\operatorname{dist}(w,w')=2\)
and \(\log K_n(v)/\log n\to1\), division by \(\log n\) gives
\(D_n/\log n\to2c_\star\) almost surely. The bound \(D_n\le2H_n\)
and the uniform integrability from Step~1 give both \(L^p\) limits for
the diameter.
\end{proof}

\begin{remark}[Storage at the leaves]\label{rem:joseph}
Joseph's results give another proof of the almost sure height limit
when the split law is non-geometric in the sense of~\cite[Section~2]{Joseph2011}
and the infimum in~\eqref{eq:transform-definition} is attained at
\(\beta_\star\). In the nested occupancy scheme~\cite{Bertoin2008}, let
\(\sigma_{n,j}=\min\{h\ge0:\max_{|u|=h}N_n(u)<j\}\) for \(j>1\).
For a fixed integer \(j\ge\beta_\star\) and \(0<\varepsilon<1\),
\cite[Theorem~1 and Proposition~1]{Joseph2011} give
\[
 \frac{\sigma_{n,j}}{\log n}\longrightarrow c_\star,
 \qquad
 \frac{\sigma_{n,n^\varepsilon}}{\log n}\longrightarrow
 (1-\varepsilon)c_\star
 \qquad\text{almost surely}.
\]

Lemma~\ref{lem:comparison} gives \(H_n\le\sigma_{n,j}+j-2\), so
\(H_n=O(\log n)\) almost surely. At depth \(H_n+1\), every vertex \(u\)
has \(K_n(u)=0\). By the same lemma,
\(\max_{|u|=H_n+1}N_n(u)\le H_n+1<n^\varepsilon\) for all sufficiently
large \(n\), almost surely. Thus \(\sigma_{n,n^\varepsilon}-1\le H_n\).
Dividing the two bounds on \(H_n\) by \(\log n\) and letting
\(\varepsilon\downarrow0\) gives \(H_n/\log n\to c_\star\) almost surely.
\end{remark}

\section{Second-order height and diameter}\label{sec:second-order}

Throughout this section, we assume that \(\Pp(M<1)>0\) and that the
infimum in~\eqref{eq:transform-definition} is attained at
\(\beta_\star\), so that~\eqref{eq:critical-constant} holds.
We write \(C\) for a finite constant that may change from line to line and is independent of \(n\), \(h\) and \(x\).

\subsection{The associated branching random walk}\label{sec:brw}

When summing over a generation, we include only vertices of positive mass. For a vertex \(u\) of positive mass let
\[
 S(u)=-\beta_\star\log\mu(u)-\frac{\beta_\star}{c_\star}\,|u| .
\]
These positions form a branching random walk in which every particle
has at least one child and, with positive probability, at least two.
By~\eqref{eq:critical-constant} it is in the boundary case, that is,
\[
 \E\sum_{|u|=1}\mathrm e^{-S(u)}
 =\mathrm e^{\beta_\star/c_\star}L(\beta_\star)=1,
 \qquad
 \E\sum_{|u|=1}S(u)\,\mathrm e^{-S(u)}
 =\kappa(\beta_\star)-\beta_\star\kappa'(\beta_\star)=0.
\]
 Since
\(\beta_\star>1\), the function \(x\mapsto x^{\beta_\star-1}(1+\log^2x)\)
is bounded on \((0,1]\). Thus for a deterministic constant \(C\),
\begin{equation}\label{eq:mallein-integrability}
 \sum_{|u|=1}\bigl(1+S(u)^2\bigr)\,\mathrm e^{-S(u)}
 \le C\sum_{i\ge1}P_i=C .
\end{equation}
In particular, \(\sum_{|u|=1}(1+|S(u)|)\mathrm e^{-S(u)}\) is bounded
by a deterministic constant. These bounds verify Madaule's moment
condition~\cite[(1.5)]{Madaule2017} for \(S\) and Mallein's moment
conditions~\cite[(1.3)--(1.4)]{Mallein2016} for \(-S\). The variance
\(\E\sum_{|u|=1}S(u)^2\mathrm e^{-S(u)}\) equals
\(\beta_\star^2\kappa''(\beta_\star)>0\) by~\eqref{eq:kappa-variance}.

For \(h\ge0\) set \(m_h=\mathrm e^{-h/c_\star}(1+h)^{-3/(2\beta_\star)}\).
The definition of \(b_n\) gives
\begin{equation}\label{eq:mass-scale}
 n\,m_h=\mathrm e^{(b_n-h)/c_\star}
 \left(\frac{\log n}{1+h}\right)^{3/(2\beta_\star)},\qquad n\ge2 .
\end{equation}
For \(2\le m\le n\), we also have
\begin{equation}\label{eq:centering-shift}
 b_m-b_n=c_\star\log\frac{m}{n}
 -\frac{3c_\star}{2\beta_\star}\log\frac{\log m}{\log n}
 \ge c_\star\log\frac{m}{n}.
\end{equation}
In particular, if \(m/n\to p\in(0,1]\) as \(n\to\infty\), then
\(b_m-b_n\to c_\star\log p\).

We use two estimates for the branching random walk, both stated in terms of the masses.
The first estimate concerns the sum \(W_h(s)=\sum_{|u|=h}\mu(u)^s\).
For every \(s>\beta_\star\), there is a constant \(C_s<\infty\) such that
\begin{equation}\label{eq:partition-tail}
 \Pp\bigl(W_h(s)\ge z\,m_h^{s}\bigr)
 \le C_s\,(1+\log z)\,z^{-\beta_\star/s},
 \qquad h\ge0,\ z\ge1.
\end{equation}
Indeed, \(m_h^{-s}W_h(s)=(1+h)^{3s/(2\beta_\star)}
\sum_{|u|=h}\mathrm e^{-(s/\beta_\star)S(u)}\).
For \(h\ge2\), Madaule's estimate~\cite[Proposition~2.1]{Madaule2017}
with exponent \(s/\beta_\star>1\) applies to the expression on the right
with \(h\) in place of \(1+h\). Since \(1+h\le2h\), this gives
\eqref{eq:partition-tail} for all sufficiently large \(z\), uniformly
in \(h\ge2\). For \(h=0,1\), the event is empty when \(z>m_h^{-s}\),
since \(W_h(s)\le1\). Increasing \(C_s\) if necessary, we obtain the bound for all \(h\ge0\) and \(z\ge1\).
For \(0<q<\beta_\star/s\), integration of~\eqref{eq:partition-tail} gives
\begin{equation}\label{eq:partition-moments}
 \sup_{h\ge0}\E\!\left[\left(\frac{W_h(s)}{m_h^s}\right)^q\right]
 \le 1+C_s q\int_1^\infty(1+\log z)\,z^{q-1-\beta_\star/s}\,dz<\infty .
\end{equation}

Second, for the maximal mass \(\mu_h^\star=\max_{|u|=h}\mu(u)\),
Mallein's theorem~\cite[Theorem~1.1]{Mallein2016} applied to \(-S\) gives
\[
 \beta_\star\log\frac{m_h}{\mu_h^\star}
 =\min_{|u|=h}S(u)-\tfrac32\log(1+h)
 =O_{\Pp}(1).
\]
Consequently,
\begin{equation}\label{eq:min-tightness}
 \frac{m_h}{\mu_h^\star}=O_{\Pp}(1).
\end{equation}
Both estimates allow infinitely many children and lattice displacements.

\subsection{Upper tail}\label{sec:upper-tail}

\begin{proof}[Proof of~\eqref{eq:second-order-upper-tail}]
Fix an integer \(j>\beta_\star\) and set \(h=\lceil b_n+x\rceil\).
For all sufficiently large \(n\), we have \(1+h\ge b_n\ge(c_\star/2)\log n\).
Since \(h\ge b_n+x\), \eqref{eq:mass-scale} gives
\(n m_h\le C_0\mathrm e^{-x/c_\star}\) for every \(x\ge0\).
Take \(t=C_0^j\mathrm e^{-jx/c_\star}\), so that \((nm_h)^j\le t\le1\)
for all sufficiently large \(x\). Since \(H_n\) is integer-valued,
taking expectations in~\eqref{eq:quenched-tail} with \(\beta=j\)
and using~\eqref{eq:partition-tail} gives
\begin{align*}
 \Pp(H_n\ge b_n+x+j-1)
 &\le \E\!\left[\min\!\left\{1,t\frac{W_h(j)}{m_h^j}\right\}\right]
 =t\int_0^{1/t}\Pp\bigl(W_h(j)>z m_h^j\bigr)\,dz\\
 &\le C t^{\beta_\star/j}\bigl(1+\log(1/t)\bigr)
 \le C(1+x)\,\mathrm e^{-\beta_\star x/c_\star}.
\end{align*}
In the integral we used the bound \(1\) for \(z\le1\) and
\(C(1+\log z)z^{-\beta_\star/j}\) for \(z>1\), with \(\beta_\star/j<1\).
Replacing \(x\) by \(x-j+1\) gives the height bound for all sufficiently
large \(n\) and \(x\). For the finitely many smaller \(n\), the event is
empty for large \(x\) since \(H_n\le n-1\), so it suffices to increase
\(C\) to cover bounded \(x\). The diameter bound follows from \(D_n\le2H_n\).
\end{proof}

\subsection{Occupation times and tightness}\label{sec:tightness}

At depth \(h=\lceil b_n\rceil\sim c_\star\log n\),
\eqref{eq:mass-scale} gives \(n m_h\asymp1\). For a vertex \(u\)
at this depth with mass comparable to \(m_h\), we thus have
\(\E[N_n(u)\mid\cE]=n\mu(u)\asymp1\).
To deduce that \(u\) is occupied from the bound \(K_n(u)\ge N_n(u)-h\)
in Lemma~\ref{lem:comparison}, we would need \(N_n(u)>h\).
We therefore use occupation times to prove the height lower bound.
For a vertex \(u\) of positive mass, let
\(\tau(u)=\inf\{n\ge1:u\in\widehat T_n\}\) be the index of the ball stored
at \(u\).

\begin{lemma}[Occupation time of a vertex of maximal mass]\label{lem:occupation}
For every vertex \(u\) of positive mass,
\begin{equation}\label{eq:mean-occupation-time}
 \E[\tau(u)\mid\cE]=\sum_{k=0}^{|u|}\frac1{\mu(u|_k)}.
\end{equation}
For each \(h\), let \(u_h^\star\) be an \(\cE\)-measurable choice of
a vertex of maximal mass at depth \(h\). Then, as \(h\to\infty\),
\begin{equation}\label{eq:occupation}
 m_h\,\E[\tau(u_h^\star)\mid\cE]=O_{\Pp}(1).
\end{equation}
\end{lemma}

\begin{proof}
\emph{Step 1. Occupation times along a path.}
The root stores ball \(1\). For \(1\le k\le |u|\), once \(u|_{k-1}\) is
occupied, \(u|_k\) stores the next ball whose ray passes through it.
Conditionally on \(\cE\), the rays after the stopping time
\(\tau(u|_{k-1})\) are independent of the past, and each hits \(u|_k\)
with probability \(\mu(u|_k)\) by~\eqref{eq:ray-law}. Hence
\(\tau(u|_k)-\tau(u|_{k-1})\) is geometric with conditional mean
\(1/\mu(u|_k)\). Summing over \(k\) gives~\eqref{eq:mean-occupation-time}.

\smallskip\noindent\emph{Step 2. A vertex of maximal mass.}
By Lemma~\ref{lem:transform}\ref{it:L} and~\eqref{eq:critical-constant},
\(L(\beta_\star+1)<L(\beta_\star)=\mathrm e^{-\beta_\star/c_\star}\).
By continuity, we can choose \(s\in(\beta_\star,2\beta_\star)\) and
\(\delta>0\) such that \(\mathrm e^{s/c_\star}L(s+1)\le\mathrm e^{-2\delta}\).
For \(0\le k\le h\), let
\[
 Z_{h,k}=m_h^{-s}\sum_{|u|=h}\frac{\mu(u)^{s+1}}{\mu(u|_k)}.
\]
Using~\eqref{eq:mean-occupation-time} and the term \(u=u_h^\star\)
in each \(Z_{h,k}\), we obtain
\begin{equation}\label{eq:occupation-bound}
 m_h\,\E[\tau(u_h^\star)\mid\cE]
 =\sum_{k=0}^h\frac{m_h}{\mu(u_h^\star|_k)}
 \le\Bigl(\frac{m_h}{\mu_h^\star}\Bigr)^{s+1}\sum_{k=0}^h Z_{h,k}.
\end{equation}
Grouping the vertices at depth \(h\) by their ancestor at depth \(k\)
and applying~\eqref{eq:subtree-moments} with \(\beta=s+1\), \(q=1\)
and \(r=h-k\), we obtain
\[
 \E[Z_{h,k}\mid\mathcal F_k]
 =m_h^{-s}\,L(s+1)^{h-k}\,W_k(s)
 \le\Bigl(\frac{1+h}{1+k}\Bigr)^{3s/(2\beta_\star)}
     \mathrm e^{-2\delta(h-k)}\,\frac{W_k(s)}{m_k^{s}}.
\]
Since \(s<2\beta_\star\), \eqref{eq:partition-moments} applies with
\(q=1/2\). By conditional Jensen's inequality and
\((1+h)/(1+k)\le1+h-k\),
\[
 \E\Bigl[\Bigl(\sum_{k=0}^{h}Z_{h,k}\Bigr)^{1/2}\Bigr]
 \le\sum_{k=0}^{h}\E\bigl[Z_{h,k}^{1/2}\bigr]
 \le C\sum_{\ell\ge0}(1+\ell)^{3s/(4\beta_\star)}\,\mathrm e^{-\delta \ell}<\infty .
\]
The last bound is uniform in \(h\), so \(\sum_{k=0}^h Z_{h,k}=O_{\Pp}(1)\).
By~\eqref{eq:min-tightness}, \(m_h/\mu_h^\star=O_{\Pp}(1)\), and
\eqref{eq:occupation-bound} gives~\eqref{eq:occupation}.
\end{proof}

\begin{proof}[Proof of~\eqref{eq:second-order-tightness}]
\leavevmode\par\nobreak\noindent
\emph{Step 1. Height.}
Given the upper-tail bound~\eqref{eq:second-order-upper-tail},
it remains to show that
\[
 \lim_{x\to\infty}\,\sup_{n\ge2}\,\Pp(H_n<b_n-x)=0.
\]
Fix \(n\ge2\). Since \(H_n\ge0\), it suffices to consider \(0\le x<b_n\).
With \(h=\lceil b_n-x\rceil\) and \(u_h^\star\) as in
Lemma~\ref{lem:occupation}, we have
\(\{H_n<b_n-x\}=\{H_n<h\}\subset\{\tau(u_h^\star)>n\}\).
Since \(h\le b_n-x+1\) and \(1+h\le b_n+2\le C\log n\),
\eqref{eq:mass-scale} gives \(n\,m_h\ge C^{-1}\mathrm e^{x/c_\star}\).
Conditioning on \(\cE\) and applying Markov's inequality gives, for \(K>0\),
\[
 \Pp(H_n<b_n-x)
 \le \E\!\left[\min\!\left\{1,\frac{\E[\tau(u_h^\star)\mid\cE]}{n}\right\}\right]
 \le C K\mathrm e^{-x/c_\star}
   +\Pp\bigl(m_h\E[\tau(u_h^\star)\mid\cE]>K\bigr).
\]
By~\eqref{eq:occupation}, the last probability tends to zero as
\(K\to\infty\), uniformly in \(h\). Taking the supremum over \(n\ge2\)
and letting first \(x\to\infty\), then \(K\to\infty\), proves the height
estimate.

\smallskip\noindent\emph{Step 2. Diameter.}
Let \(w,w'\) be the two children selected in Step~2 of the proof of
Theorem~\ref{thm:height}.
For each \(v\in\{w,w'\}\), Lemma~\ref{lem:comparison} gives
\(K_n(v)/n\to\mu(v)>0\) almost surely. Hence \(K_n(v)\ge2\)
eventually, and~\eqref{eq:centering-shift} gives
\(b_{K_n(v)}-b_n\to c_\star\log\mu(v)\).
By Lemma~\ref{lem:branching}, for every \(x\ge0\),
\[
 \Pp\!\left(K_n(v)\ge2,\,
   \bigl|H(\widehat T^{\,v}_{K_n(v)})-b_{K_n(v)}\bigr|>x\right)
 \le \sup_{m\ge2}\Pp(|H_m-b_m|>x).
\]
The right-hand side tends to zero as \(x\to\infty\) by Step~1.
Thus both subtree heights are \(b_n+O_{\Pp}(1)\), and
\eqref{eq:two-subtree-diameter} gives \(D_n=2b_n+O_{\Pp}(1)\).
\end{proof}

\subsection{Lower tails and moments}\label{sec:lower-tail}

\begin{proof}[Proof of Theorem~\ref{thm:lower-tail}]
\leavevmode\par\nobreak\noindent
\emph{Step 1. A lower bound from the root split.}
Write \(P=P(\varnothing)\) for the root split and set \(Y=-\log M\).
We first show that, for some \(A>0\) and all sufficiently large \(n\),
\begin{equation}\label{eq:root-deficit}
 \Pp(b_n-H_n>x)\ge\tfrac12\Pp\bigl(Y\ge A(1+x)\bigr),
 \qquad 0\le x\le b_n/2 .
\end{equation}
To bound \(\Pp(H_n\ge b_n-x\mid P)\) from above, fix an integer
\(j>\beta_\star\) and choose \(1<a<\beta_\star\).
For \(h\ge1\), decompose the sum at the root as
\(W_h(j)=\sum_{i:P_i>0}\sum_{|w|=h-1}\mu(iw)^j\).
We use this decomposition in~\eqref{eq:quenched-tail} with \(\beta=j\).
Since \(0<a/j<1\), we have
\(\min\{1,\sum_i t_i\}\le(\sum_i t_i)^{a/j}\le\sum_i t_i^{a/j}\)
for \(t_i\ge0\). Taking the conditional expectation given \(P\) and
applying~\eqref{eq:subtree-moments} with \(k=1\) and \(r=h-1\) gives
\[
 \Pp(H_n\ge h+j-1\mid P)
 \le n^a\sum_iP_i^a\,\E[W_{h-1}(j)^{a/j}]
 \le C(nm_{h-1})^a M^{a-1}.
\]
The last inequality uses~\eqref{eq:partition-moments} with \(s=j\)
and exponent \(a/j<\beta_\star/j\), together with the bound
\(\sum_iP_i^a\le M^{a-1}\) from Lemma~\ref{lem:transform}\ref{it:L}.
For \(0\le x\le b_n/2\), set \(h=\lceil b_n-x\rceil-j+1\). Then
\(b_n-(h-1)\le x+j\) and, for all sufficiently large \(n\),
\(h\ge(c_\star/4)\log n\) uniformly in \(x\). Applying~\eqref{eq:mass-scale}
at depth \(h-1\) gives \(nm_{h-1}\le C\mathrm e^{x/c_\star}\).
Thus \(\Pp(H_n\ge b_n-x\mid P)\le C\exp(ax/c_\star-(a-1)Y)\).
For \(A\) sufficiently large, this bound gives
\(\Pp(b_n-H_n>x\mid P)\ge\tfrac12\) whenever \(Y\ge A(1+x)\).
Averaging over \(P\) yields~\eqref{eq:root-deficit}.

Integrating~\eqref{eq:root-deficit} against \(px^{p-1}\,dx\) gives
\[
 \E[(b_n-H_n)_+^p]\ge
 \tfrac12\,\E\!\left[\bigl((Y/A-1)_+\wedge(b_n/2)\bigr)^p\right].
\]
Since \(b_n\sim c_\star\log n\), for all large \(n\) we have
\[
 Y\wedge\log n
 \le A+C\bigl((Y/A-1)_+\wedge(b_n/2)\bigr).
\]
Taking \(p\)th moments in the last inequality and using the preceding
estimate, we obtain
\(1+\E[(b_n-H_n)_+^p]\ge C^{-1}\bigl(1+\E[(Y\wedge\log n)^p]\bigr)\).
Since \(D_n/2\le H_n\), the same lower bound holds for the diameter
term in~\eqref{eq:truncated-moments}.

We next show that either bound in~\eqref{eq:lower-tail}
implies~\eqref{eq:moment-assumption}. If the height bound in~\eqref{eq:lower-tail} holds, then for each
\(x\ge0\) we may choose \(n\) with \(b_n\ge2x\) in~\eqref{eq:root-deficit}
to get \(\Pp(Y\ge A(1+x))\le C\mathrm e^{-\delta x}\).
Hence \(\E[M^{-\eta}]=\E[\mathrm e^{\eta Y}]<\infty\) for
\(0<\eta<\delta/A\). The diameter bound implies the height bound since
\(D_n\le2H_n\).

\smallskip\noindent\emph{Step 2. A recursion for the diameter.}
Let \(w_1,w_2\) be the children corresponding to the two largest
coordinates \(M\ge M'\) of the root split, breaking ties by index.
Write \(K_i=K_n(w_i)\) and let \(D^{(i)}\) be the diameter of the
subtree at \(w_i\), with diameter zero for an empty subtree. In this
proof set \(b_0=b_1=D_0=0\). By Lemma~\ref{lem:branching}, given
\(P\) and \(K_1=k_1,K_2=k_2\), the diameters \(D^{(1)},D^{(2)}\)
are independent with the laws of \(D_{k_1},D_{k_2}\), respectively.

Since \(\Pp(M'>0)=\Pp(M<1)>0\), choose \(\varepsilon>0\) with
\(\Pp(M'>\varepsilon)>0\). By~\eqref{eq:second-order-tightness}, there is
\(B>0\) such that \(\inf_{m\ge2}\Pp(D_m/2\ge b_m-B)\ge1/2\). Set
\(G_n=\{K_2\ge\varepsilon n,\,D^{(2)}/2\ge b_{K_2}-B\}\).
Since \(K_2/n\to M'\) almost surely by Lemma~\ref{lem:comparison},
the conditional law of \(D^{(2)}\) gives
\[
 \liminf_{n\to\infty}\Pp(G_n)
 \ge\tfrac12\liminf_{n\to\infty}\Pp(K_2\ge\varepsilon n)
 \ge\tfrac12\Pp(M'>\varepsilon)>0 .
\]
Fix \(n_0\) and \(\varrho\in(0,1)\) such that \(\varepsilon n_0\ge2\)
and \(\Pp(G_n)\ge\varrho\) for \(n\ge n_0\).

Both subtree diameters are at most \(D_n\). On \(G_n\), for \(n\ge n_0\),
\eqref{eq:centering-shift} gives
\(b_n-D_n/2\le b_n-b_{K_2}+B\le B+c_\star\log(1/\varepsilon)\).
If \(K_1\ge2\), \eqref{eq:centering-shift} with \(m=K_1\) gives
\(b_n-D_n/2\le c_\star\log(n/K_1)+b_{K_1}-D^{(1)}/2\).
If \(K_1\le1\), then \(b_{K_1}=D^{(1)}=0\) and
\(b_n\le c_\star\log n+C\). Set \(r_n=n/(K_1+1)\). Since
\(K_1+1\le2(K_1\vee1)\), there is a constant \(C_0>0\) such that
\begin{equation}\label{eq:diameter-recursion}
 b_n-D_n/2
 \le C_0+\mathbf1_{G_n^c}\bigl(c_\star\log r_n+b_{K_1}-D^{(1)}/2\bigr),
 \qquad n\ge n_0 .
\end{equation}

\smallskip\noindent\emph{Step 3. Truncated moments and exponential lower tails.}
We first estimate \(\log r_n\) in~\eqref{eq:diameter-recursion} using
the conditional law of \(K_1\). Ball \(1\) is stored at the root. Given \(P\), the remaining \(n-1\) balls
enter the subtree at \(w_1\) independently with probability \(M\), so
\(K_1\sim\Bin(n-1,M)\). Using
\(\frac{n}{k+1}\binom{n-1}{k}=\binom n{k+1}\), we obtain
\begin{equation}\label{eq:root-count-reciprocal}
 \E[M r_n\mid P]=1-(1-M)^n\le1 .
\end{equation}
Since \(r_n\le n\), truncating \(Y\) at \(\log n\) leaves the event
\(\{\log r_n>Y+x\}\) unchanged for \(x\ge0\). Hence
\eqref{eq:root-count-reciprocal} and Markov's inequality give
\[
 \Pp\bigl(\log r_n>(Y\wedge\log n)+x\mid P\bigr)
 =\Pp(Mr_n>\mathrm e^x\mid P)\le\mathrm e^{-x}.
\]
Fix \(p>0\). Integration against \(px^{p-1}\,dx\) gives
\(\E[(\log r_n-(Y\wedge\log n))_+^p]\le\Gamma(p+1)\), hence
\(\E[(\log r_n)^p]\le C_p(1+\E[(Y\wedge\log n)^p])\).

To take expectations in~\eqref{eq:diameter-recursion}, let \(f\) be a
nonnegative function and let \(Z\ge0\) be measurable with respect to
\(\sigma(P,K_1,K_2,D^{(2)})\). The conditional independence in Step~2 gives
\begin{equation}\label{eq:root-subtree-expectation}
 \E[Z f(K_1,D^{(1)})]
 =\sum_{m=0}^{n-1}\E[Z\mathbf1_{\{K_1=m\}}]\,\E[f(m,D_m)]
 \le \E[Z]\max_{0\le m<n}\E[f(m,D_m)] .
\end{equation}
Choose \(1<\lambda<(1-\varrho)^{-1}\). For \(u\ge0\) and \(v\in\mathbb R\),
subadditivity when \(p\le1\) and convexity when \(p>1\) give
\((u+v)_+^p\le C_{p,\lambda}u^p+\lambda(v_+)^p\).
Apply this bound in~\eqref{eq:diameter-recursion} and
take \(Z=\mathbf1_{G_n^c}\), \(f(m,d)=(b_m-d/2)_+^p\)
in~\eqref{eq:root-subtree-expectation}. Since \(\Pp(G_n^c)\le1-\varrho\),
\[
 \E[(b_n-D_n/2)_+^p]
 \le C\bigl(1+\E[(Y\wedge\log n)^p]\bigr)
     +\lambda(1-\varrho)\max_{0\le m<n}\E[(b_m-D_m/2)_+^p].
\]
Divide by \(1+\E[(Y\wedge\log n)^p]\), which is non-decreasing in \(n\).
The terms with \(m=0,1\) vanish. Since \(\lambda(1-\varrho)<1\),
induction on \(n\) gives
\(\E[(b_n-D_n/2)_+^p]\le C\bigl(1+\E[(Y\wedge\log n)^p]\bigr)\)
for all \(n\ge2\). Since \(D_n/2\le H_n\), the same upper bound holds
for the height. Together with the lower bounds from Step~1, this
gives~\eqref{eq:truncated-moments} after adjusting the constants for
the finitely many small \(n\). By monotone convergence,
\(\E[(Y\wedge\log n)^p]\uparrow\E[Y^p]\). Since
\eqref{eq:second-order-upper-tail} gives uniform \(p\)th moment bounds
for \((H_n-b_n)_+\) and \((D_n/2-b_n)_+\), we obtain
\eqref{eq:centred-moments}.

To deduce~\eqref{eq:lower-tail} from~\eqref{eq:moment-assumption}, fix
\(\eta>0\) with \(\E[M^{-\eta}]<\infty\). For \(0<tc_\star\le\min\{1,\eta\}\),
conditional Jensen's inequality and~\eqref{eq:root-count-reciprocal}
give \(\E[r_n^{tc_\star}\mid P]\le M^{-tc_\star}\). Since \(r_n\ge1\) and \(\Pp(G_n)\ge\varrho\) for \(n\ge n_0\),
\[
 \E[r_n^{tc_\star}\mathbf1_{G_n^c}]
 \le \E[r_n^{tc_\star}]-\Pp(G_n)
 \le \E[M^{-tc_\star}]-\varrho .
\]
Applying \(x\mapsto\mathrm e^{tx}\) to~\eqref{eq:diameter-recursion} and using
\eqref{eq:root-subtree-expectation} with \(Z=r_n^{tc_\star}\mathbf1_{G_n^c}\)
and \(f(m,d)=\mathrm e^{t(b_m-d/2)}\) gives
\[
 \E[\mathrm e^{t(b_n-D_n/2)}]
 \le \mathrm e^{tC_0}
  +\mathrm e^{tC_0}\bigl(\E[M^{-tc_\star}]-\varrho\bigr)
     \max_{0\le m<n}\E[\mathrm e^{t(b_m-D_m/2)}].
\]
By dominated convergence, the coefficient of the maximum tends to
\(1-\varrho<1\) as \(t\downarrow0\). For sufficiently small fixed \(t>0\),
induction gives \(\sup_{n\ge2}\E[\mathrm e^{t(b_n-D_n/2)}]<\infty\).
Markov's inequality gives the diameter bound in~\eqref{eq:lower-tail},
and \(D_n\le2H_n\) gives the height bound.
\end{proof}

\begin{proof}[Proof of~\eqref{eq:variance-order}]
\leavevmode\par\nobreak\noindent
\emph{Step 1. Reduction to a diameter estimate.}
For fixed \(n\ge2\), let \(X\) denote either \(H_n-b_n\) or \(D_n/2-b_n\),
and write \(X_-=(-X)_+\). By~\eqref{eq:second-order-tightness}, choose
\(B>0\) such that \(\Pp(|X|\le B)\ge1/2\) for both choices and all \(n\ge2\).
The upper-tail bound~\eqref{eq:second-order-upper-tail} gives
\(\Var(X)\le\E[X^2]\le C+\E[X_-^2]\).
For an independent copy \(X'\) of \(X\), symmetry and independence give
\begin{align*}
 \Var(X)=\E[(X'-X)_+^2]
 &\ge \Pp(X'\ge B+1)\,\E[(B+1-X)_+^2]\\
 &\ge \Pp(X\ge B+1)\bigl(\tfrac12+\E[X_-^2]\bigr),
\end{align*}
where \((B+1-x)_+^2\ge\mathbf1_{\{x\le B\}}+x_-^2\).
Since \(X\ge D_n/2-b_n\) in both cases, these bounds and
\eqref{eq:truncated-moments} with \(p=2\) reduce the proof to showing
that \(\liminf_{n\to\infty}\Pp(D_n/2\ge b_n+B+1)>0\).

\smallskip\noindent\emph{Step 2. Two subtrees in distinct root branches.}
Let \(M'\) be the second largest coordinate of the root split and choose
\(\varepsilon>0\) with \(\Pp(M'\ge\varepsilon)>0\).
By Lemma~\ref{lem:transform}\ref{it:attained}, fix
\(t\in(\mathrm e^{-1/c_\star},p_{\max})\). Since \(1+c_\star\log t>0\),
choose an integer \(r\ge1\) such that \(r+1+c_\star\log(\varepsilon t^r)>2B+1\).
From each of the two children with the largest root split fractions, take \(r\)
further steps, always choosing a child with the largest split fraction.
Break ties by index and let \(w_1,w_2\) be the endpoints.
Each vertex on the two paths is chosen using only splits at its strict
ancestors. Thus, given the root split, the \(2r\) split fractions chosen below it
are independent and each has the law of \(M=\max_{i\ge1}P_i\). In particular, the event
\(E_r=\{\min_{i=1,2}\mu(w_i)\ge\varepsilon t^r\}\) satisfies
\[
 \Pp(E_r)\ge\Pp(M'\ge\varepsilon)\,\Pp(M\ge t)^{2r}>0 .
\]

On \(E_r\), Lemma~\ref{lem:comparison} gives
\(K_n(w_i)/n\to\mu(w_i)\ge\varepsilon t^r>0\) almost surely for \(i=1,2\).
Thus \(\min_{i=1,2}K_n(w_i)\ge2\) for all sufficiently large \(n\), and
\eqref{eq:centering-shift} gives
\(r+1+b_{K_n(w_i)}-b_n\to r+1+c_\star\log\mu(w_i)>2B+1\).
The event \(E_r\) and the two subtree sizes \(K_n(w_1)\) and \(K_n(w_2)\) are measurable with respect
to the \(\sigma\)-field \(\mathcal G\) in Lemma~\ref{lem:branching}.
Almost surely on \(E_r\), for all sufficiently large \(n\),
\eqref{eq:two-subtree-diameter} and \(\operatorname{dist}(w_1,w_2)=2(r+1)\)
give \(D_n/2\ge b_n+B+1\) whenever the subtree at \(w_i\) has height
at least \(b_{K_n(w_i)}-B\) for \(i=1,2\). Conditional independence
therefore gives
\[
 \Pp(D_n/2\ge b_n+B+1\mid\mathcal G)
 \ge \left(\inf_{m\ge2}\Pp(H_m\ge b_m-B)\right)^2
 \ge\tfrac14 ,
\]
where the last inequality follows from the choice of \(B\) in Step~1.
Fatou's lemma gives
\(\liminf_{n\to\infty}\Pp(D_n/2\ge b_n+B+1)\ge\Pp(E_r)/4>0\).
\end{proof}

\begin{remark}[Two examples with unbounded branching]\label{rem:mean-height}
The centered height and diameter can be tight even when their means tend
to \(-\infty\) and their variances tend to infinity. For example, take a
\(\GEM(1)\) split with probability \(1/2\) and otherwise divide the mass
equally among \(\lceil\mathrm e^V\rceil\) children, where
\(\Pp(V>y)=(1+y)^{-1/2}\) for \(y\ge0\). Here \(p_{\max}=1\) and
\(\Pp(M=1)=0\), so~\eqref{eq:attainment} and Theorem~\ref{thm:second-order}
give tightness of \(H_n-b_n\) and \(D_n-2b_n\).
Write \(Y=-\log M\). On the equal-split event, \(Y=\log\lceil\mathrm e^V\rceil\), so \(V\le Y\le V+\log 2\).
On the \(\GEM(1)\) event, \(Y\le-\log P_1\) with \(P_1\) uniform on
\((0,1)\). It follows that \(\Pp(Y>y)\sim1/(2\sqrt y)\), so for \(p>1/2\),
\[
 \E[(Y\wedge t)^p]
 =p\int_0^t y^{p-1}\Pp(Y>y)\,dy
 \sim \frac{p}{2p-1}\,t^{p-1/2}
 \qquad(t\to\infty).
\]
We apply this formula with \(t=\log n\) for \(p=1\) and \(p=2\).
The mean estimates following Theorem~\ref{thm:lower-tail} and the
variance estimate~\eqref{eq:variance-order} then give
\[
 b_n-\E H_n \asymp 2b_n-\E D_n \asymp \sqrt{\log n},
 \qquad
 \Var(H_n) \asymp \Var(D_n) \asymp (\log n)^{3/2}.
\]

To compare the moment and exponential lower-tail conditions in
Theorem~\ref{thm:lower-tail}, now take \(\Pp(V>y)=\mathrm e^{-\sqrt y}\)
for \(y\ge0\). The same bounds on \(Y\) give \(\E[Y^p]<\infty\) for
every \(p>0\).
On the equal-split event \(Y\ge V\), so
\(\E[\mathrm e^{\eta Y}]\ge\tfrac12\E[\mathrm e^{\eta V}]=\infty\)
for every \(\eta>0\). By Theorem~\ref{thm:lower-tail}, the absolute moments
in~\eqref{eq:centred-moments} are uniformly bounded for every \(p>0\),
while neither lower tail has a uniform exponential bound.
\end{remark}

\printbibliography
\end{document}